\documentclass[11pt, amssymb]{amsart}
\usepackage{amsmath, amsthm, amsfonts, verbatim, amssymb}
\usepackage{eucal}

\DeclareFontFamily{OT1}{rsfs}{}
\usepackage[all]{xy}

\DeclareFontShape{OT1}{rsfs}{n}{it}{<-> rsfs10}{}
\DeclareMathAlphabet{\mathscr}{OT1}{rsfs}{n}{it}
\def\sG{{\mathscr G}}
\def\fm{\mathfrak m}

\def\osG{\overline{\mathscr{G}}}

\def\osH{\overline{\mathscr{H}}}
\def\osP{\overline{\mathscr{P}}}
\def\osT{\overline{\mathscr{T}}}
\def\osQ{\overline{\mathscr{Q}}}
\def\orG{\overline{G}^{\,\rm pred}}
\def\orH{\overline{\mathscr{H}}^{\,\rm pred}}
\def\opredH{\overline{H}^{\,\rm pred}}
\def\cB{\mathscr B}

\def\sH{\mathscr H}

\def\sR{{\mathscr R}}

\def\hcO{\widehat{\mathcal{O}}}
\def\hK{\widehat{K}}
\def\hh{\widehat{h}}

\def\bR{{\mathbb R}}
\def\bZ{{\mathbb Z}}
\def\bQ{{\mathbb Q}}

\def\cB{{\mathcal B}}

\def\sS{{\mathscr S}}
\def\cG{{\mathcal G}}
\def\cZ{{\mathcal Z}}
\def\sH{{\mathscr H}}
\def\cS{{\mathcal S}}
\def\sT{{\mathscr T}}

\def\cC{{\mathcal C}}
\def\cX{{\mathcal X}}
\def\sX{{\mathscr X}}

\def\ni{\noindent}
\def\bQ{{\mathbb Q}}

\def\osH{\overline{\mathscr{H}}}
\def\osK{\overline{\mathscr{K}}}

\def\fK{\mathfrak{K}}
\def\osS{{\overline{\mathscr{S}}}}
\def\osT{{\overline{\mathscr{T}}}}

\def\fo{\mathfrak{o}}
\def\fp{\mathfrak{p}}
\def\fT{\mathfrak{T}}
\def\sB{\mathscr{B}}
\def\fD{\mathfrak{D}}

\def\fG{\mathfrak{G}}
\def\fL{\mathfrak{L}}

\def\fZ{\mathfrak{Z}}
\def\cO{\mathcal{O}}
\def\cA{\mathscr{A}}
\def\cF{\mathcal{F}}
\def\cC{\mathcal{C}}
\def\sU{\mathscr{U}}

\def\sX{\mathscr{X}}
\def\sD{\mathscr{D}}
\def\cP{\mathcal{P}}
\def\cQ{\mathcal{Q}}

\def\cS{\mathcal{S}}
\def\cT{\mathcal{T}}
\def\cU{\mathcal{U}}

\def\cH{\mathcal{H}}
\def\cM{\mathcal{M}}
\def\cN{\mathcal{N}}
\def\cR{\mathcal{R}}
\def\ocH{\overline{\mathcal{H}}}

\def\ocS{\overline{\mathcal{S}}}

\begin{document}
\title[Finite group actions  and tamely-ramified descent]
{Finite group actions on reductive groups and buildings and tamely-ramified descent in Bruhat-Tits theory}
\maketitle
{\centerline {{By} {\sc {Gopal Prasad}}}
\vskip4mm

\centerline{\it Dedicated to Guy Rousseau}
\vskip4mm

\begin{abstract}  Let  $K$  be a discretely valued field with 
Henselian valuation ring and separably closed (but not necessarily perfect) residue field of characteristic $p$, $H$ a connected reductive $K$-group, and $\Theta$ a finite group of automorphisms of  $H$. We assume that  $p$ does not divide the order of $\Theta$ and  
Bruhat-Tits theory is available for $H$ over $K$ with $\cB(H/K)$ the Bruhat-Tits building of $H(K)$. We will show that then Bruhat-Tits theory is also available  for 
$G := (H^{\Theta})^{\circ}$ and $\cB(H/K)^{\Theta}$ is the Bruhat-Tits building of $G(K)$. (In case the residue field of $K$ is perfect, this result was proved in [PY1] by a different method.) 
As a consequence of  this result, we obtain that if Bruhat-Tits theory is available for a connected reductive $K$-group $G$ over a finite tamely-ramified extension $L$ of $K$, then it is also available for $G$ over $K$ 
and $\cB(G/K) = \cB(G/L)^{{\rm{Gal}}(L/K)}$. Using this, we prove that if $G$ is quasi-split over $L$, then it is already quasi-split over $K$.
\end{abstract}  
\vskip4mm

\ni{\bf Introduction.} This paper is a sequel to our recent paper [P2]. We will assume familiarity with that paper and  will freely 
use notions and results given in it. 
\vskip1mm

Let $\cO$ be a discretely valued Henselian local ring with valuation $\omega$. Let $\fm$ be the maximal ideal of $\cO$ and $K$ the field of fractions of $\cO$. We will assume throughout that the residue field $\kappa$ of $\cO$ is separably closed. Let $\hcO$ denote the completion of $\cO$ with respect to the valuation $\omega$ and $\hK$ the completion of $K$. For any $\cO$-scheme $\sX$, $\sX(\cO)$ and $\sX(\hcO)$ will always be assumed to carry the Hausdorff-topology induced from the metric-space topology on $\cO$ and $\hcO$ respectively. It is known that if $\sX$ is smooth, then $\sX(\cO)$ is dense in $\sX(\hcO)$, [GGM, Prop.\,3.5.2]. Similarly, for any $K$-variety $\cX$, $\cX(K)$ and $\cX(\hK)$ will be assumed to carry the Hausdorff-topology induced from the metric-space topology on $K$ and $\hK$ respectively.  In case $\cX$ is a  smooth $K$-variety, $\cX(K)$ is dense in $\cX(\hK)$, [GGM, Prop.\,3.5.2]. 
\vskip1mm

Throughout this paper $H$ will denote a connected reductive $K$-group and $\sD(H)$ its derived subgroup. In this introduction,  and beginning with \S2 everywhere, we will assume that Bruhat-Tits theory is available for $H$ over $K$ [P2, 1.9,\,1.10].  It is known (cf.\,[P2,\:3.11, 1.11]) that then Bruhat-Tits theory  is also available over $K$ for the centralizer of any $K$-split torus in $H$ and for the derived subgroup of such centralizers.  Thus there is an affine building called the Bruhat-Tits building of $H(K)$,  that is a polysimplicial complex given with a metric,  and $H(K)$ acts on it by polysimplicial isometries.  This building is also the Bruhat-Tits building of $\sD(H)(K)$ and we will denote it by $\cB(\sD(H)/K)$.  
\vskip.5mm

Let $\fZ$ be the maximal $K$-split torus  in the center of $H$. Let $V(\fZ) = \bR\otimes_{\bZ}{\rm{Hom}}_K({\rm{GL}}_1,\fZ)$. Then there is a natural action of $H(K)$ on this Euclidean space by translations, with $\sD(H)(K)$, as well as every bounded subgroup of $H(K)$,  acting trivially\,(cf.\,[P2, 1.9]).  The {\it enlarged} Bruhat-Tits building $\cB(H/K)$ of $H(K)$ is the product $V(\fZ)\times \cB(\sD(H)/K)$.  The apartments of this building, as well as that of $\cB(\sD(H)/K)$, are in bijective correspondence with maximal $K$-split tori of $H$.  Given a maximal $K$-split torus $T$ of $H$, the corresponding apartment of $\cB(H/K)$ is an affine space under $V(T) := \bR\otimes_{\bZ}{\rm{Hom}}_K({\rm{GL}}_1,T)$. 
\vskip.5mm

Given a nonempty bounded subset $\Omega$ of an apartment of $\cB(H/K)$, or of $\cB(\sD(H)/K)$, there is a smooth affine $\cO$-group scheme $\sH_{\Omega}$ with generic fiber $H$, associated to $\Omega$, such that $\sH_{\Omega}(\cO)$ is the subgroup $H(K)^{\Omega}$ of $H(K)$ consisting of elements that fix $\Omega$ and $V(\fZ)$ pointwise.  The neutral component $\sH^{\circ}_{\Omega}$ of $\sH_{\Omega}$ is an open affine $\cO$-subgroup scheme of  the latter; it is by definition the union of the generic fiber $H$  of $\sH_{\Omega}$ and the identity component of its special fiber.  The group scheme $\sH^{\circ}_{\Omega}$ is called the Bruhat-Tits  group scheme associated to $\Omega$.  It has  the properties described in [P2,\,1.9\,-1.10]. The special fiber of $\sH^{\circ}_{\Omega}$ will be denoted by $\osH^{\circ}_{\Omega}$. 
\vskip1mm

Let $\Theta$ be a finite group of automorphisms of $H$. We assume that the order of $\Theta$ is not divisible by the characteristic of the residue field $\kappa$. Let $G = (H^{\Theta})^{\circ}$.  This group is also reductive, see [Ri, Prop.\,10.1.5] or [PY1, Thm.\,2.1]. The goal of this paper is to show that Bruhat-Tits theory is available for $G$ over $K$, and the {\it enlarged} Bruhat-Tits building of $G(K)$ can be identified with  the subspace  $\cB(H/K)^{\Theta}$ of $\cB(H/K)$ consisting of points fixed under $\Theta$ (see \S3).   These results have been inspired by the main theorem of [PY1], which implies that if the residue field $\kappa$ is algebraically closed (then every reductive $K$-group is quasi-split [P2,\,1.7], so Bruhat-Tits theory is available for any such group over $K$), the enlarged Bruhat-Tits building of $G(K)$ is indeed $\cB(H/K)^{\Theta}$. 
 \vskip.5mm

In \S4, we will use the above results to obtain ``tamely-ramified descent'':  (1) We will show that if a connected reductive $K$-group $G$ is quasi-split over a finite tamely-ramified extension $L$ of  $K$, then  it is quasi-split over $K$ (Theorem 4.4); this result  has been proved by Philippe Gille in [Gi] by an entirely different method. (2) The enlarged Bruhat-Tits building $\cB(G/K)$ of $G(K)$ can be identified with the subspace of points of the enlarged Bruhat-Tits building of $G(L)$ that are fixed under the action of the Galois group ${\rm{Gal}}(L/K)$. This latter result was proved by Guy Rousseau in his unpublished thesis\,[Rou, Prop.\,5.1.1]. It is a pleasure to dedicate this paper to him for his important contributions to Bruhat-Tits theory. 
\vskip2mm

\ni{\it Acknowledgements.} I thank Brian Conrad, Bas Edixhoven and Philippe Gille for their helpful comments. I thank the referee for carefully reading the paper 
and for her/his detailed comments and suggestions which helped to improve the exposition.   
I was partially supported by NSF-grant DMS-1401380.
\vskip2mm

For a $K$-split torus $S$, let ${\rm X}_{\ast}(S)= {\rm{Hom}}({\rm{GL}}_1, S)$ and $V(S) : = \bR\otimes_{\bZ}{\rm X}_{\ast}(S)$. Then for a maximal $K$-split torus $T$ of $H$, the apartment $A(T)$ of $\cB(H/K)$ corresponding to $T$ is an affine space under $V(T)$.   
\vskip4mm

\ni{\bf 1. Passage to completion} 
\vskip2mm

We begin by proving the following well-known result.
\vskip2mm

\ni{\bf Proposition 1.1.} {\em $K$-{\rm{rank}}\,$H$ = $\hK$-{\rm{rank}}\,$H$.} 
\vskip2mm

\ni{\it Proof.} Let $T$ be a maximal $K$-split torus of $H$ and $Z$ be its centralizer in $H$. Let $Z_{\rm a}$ be the maximal $K$-anisotropic connected normal subgroup of $Z$. Then
\vskip1mm

$\hK$-{\rm{rank}}\,$H = \hK$-{\rm{rank}}\,$Z = \dim(T)+$ $\hK$-{\rm{rank}}\,$Z_{\rm a}=$ $K$-{\rm{rank}}\,$H+$$\hK$-{\rm{rank}}\,$Z_{\rm a}$ .   
\vskip1mm

\ni So to prove the proposition, it suffices to show that $Z_{\rm a}$ is anisotropic over $\hK$. But according to Theorem 1.1 of [P2], $Z_{\rm a}$ is anisotropic over $\hK$ if and only if $Z_{\rm a}(\hK)$ is bounded.  The same theorem implies that $Z_{\rm a}(K)$ is bounded. As $Z_{\rm a}(K)$ is dense in $Z_{\rm a}(\hK)$, we see that $Z_{\rm a}(\hK)$ is bounded. \hfill$\Box$

\vskip1mm

\ni{\bf Proposition 1.2.} {\em Bruhat-Tits theory for $H$ is available over $K$ if and only if it is available over $\hK$. Moreover, if Bruhat-Tits theory for $H$ is available over $K$, then the enlarged Bruhat-Tits buildings of $H(K)$ and $H(\hK)$ are equal.}
\vskip 1mm

It was shown by Guy Rousseau in his thesis that the enlarged Bruhat-Tits buildings of $H(K)$ and $H(\hK)$ coincide\,[Rou, Prop.\,2.3.5]. Moreover, every apartment in the building of $H(K)$ is also an apartment in the building of $H(\hK)$; however, the latter may have many more apartments.

\vskip2mm

\ni{\it Proof.}  We assume first that Bruhat-Tits theory is available for $H$ over $K$ and let $\cB(\sD(H)/K)$ and $\cB(H/K)$ respectively denote the Bruhat-Tits building and the enlarged Bruhat-Tits building of $H(K)$. We begin by showing that  the actions of $H(K)$ on $\cB(\sD(H)/K)$ and $\cB(H/K)$ extend to actions of $H(\hK)$ by isometries. For this purpose, we recall that $H(K)$ is dense in $H(\hK)$ and the isotropy at any point of the two buildings is a bounded open subgroup of $H(K)$.  Now let $\{h_i\}$ be a sequence in $H(K)$ which converges to a point $\hh\in H(\hK)$, then given any open subgroup of $H(K)$, for all large $i$ and $j$, $h_i^{-1}h_j$ lies in this open subgroup. Thus for any point $x$ of $\cB(\sD(H)/K)$ or $\cB(H/K)$, the sequence $h_i\cdot x$ is eventually constant, i.e., there exists a positive integer $n$ such that  $h_i\cdot x = h_n\cdot x$ for all $i\geqslant n$. We define $\hh\cdot x = h_n\cdot x$. This gives well-defined  actions of $H(\hK)$ on $\cB(\sD(H)/K)$ and $\cB(H/K)$ by isometries. 
\vskip.5mm

For a nonempty bounded subset $\Omega$ of an apartment of $\cB(H/K)$, or of $\cB(\sD(H)/K)$, let $\sH_{\Omega}$ and $\sH^{\circ}_{\Omega}$ be the smooth affine $\cO$-group schemes as in the Introduction. Then as $\sH_{\Omega}(\hcO)$ is a closed and open subgroup of $H(\hK)$ containing $\sH_{\Omega}(\cO)$ as a dense subgroup, we see that  $\sH_{\Omega}(\hcO)$ equals the subgroup $H(\hK)^{\Omega}$ of $H(\hK)$ consisting of elements that fix  $\Omega$, as well as $V(\fZ)$, pointwise.   
\vskip.5mm

Let $T$ be a maximal $K$-split torus of $H$, then by Proposition 1.1, $T_{\hK}$ is a maximal $\hK$-split torus of $H_{\hK}$. Let $A$ be the apartment of $\cB(\sD(H)/K)$, or of $\cB(H/K)$,  corresponding to $T$. Then every maximal $\hK$-split torus of $H_{\hK}$ is of the form $\hh T_{\hK}\hh^{-1}$ for an $\hh\in H(\hK)$, and we define the corresponding apartment to be $\hh\cdot A$. We now declare $\cB(\sD(H)/K)$ (resp.\,$\cB(H/K)$) to be the Bruhat-Tits building (resp.\:the enlarged Bruhat-Tits building) of $H(\hK)$ with these apartments. 
\vskip.5mm 

Let $A$ be an apartment of $\cB(\sD(H)/K)$, or of $\cB(H/K)$,   corresponding to a maximal $K$-split torus $T$ of $H$ and $\hh\in H(\hK)$. Given a nonempty bounded subset $\widehat\Omega$ of $\widehat{A}:=\hh\cdot A$, the subset $\Omega:= \hh^{-1}\cdot{\widehat\Omega}$ is contained in $A$. The closed and open subgroup $\hh H(\hK)^{\Omega}\hh^{-1} = \hh\sH_{\Omega}(\hcO)\hh^{-1}$ of $H(\hK)$ is the subgroup $H(\hK)^{\widehat\Omega}$ consisting  of elements that fix $\widehat\Omega$, as well as $V(\fZ)$, pointwise.  Now as $H(K)$ is dense in $H(\hK)$ and $H(\hK)^{\widehat\Omega}$ is an open subgroup, $H(\hK) = H(\hK)^{\widehat\Omega}\cdot H(K)$, so $\hh = h'\cdot h$, with $h'\in H(\hK)^{\widehat\Omega}$ and $h\in H(K)$. Thus the apartment $\widehat{A} =\hh\cdot A=  h'\cdot hA$, and $hA$ is an apartment of the Bruhat-Tits building of $H(K)$ (or of the enlarged building of $H(K)$). As $h'\in H(\hK)^{\widehat\Omega}$, the apartment $hA$ contains $\widehat\Omega$. This shows that any bounded subset $\widehat\Omega$ of an apartment of  the Bruhat-Tits building of $H(\hK)$ (or of the enlarged building of $H(\hK)$) is contained in an apartment of the Bruhat-Tits building of $H(K)$ (or of the enlarged building of $H(K)$). We define the $\hcO$-group schemes $\sH_{\widehat\Omega}$ and $\sH^{\circ}_{\widehat\Omega}$ associated to $\widehat\Omega$ to be the group schemes obtained from the corresponding $\cO$-group schemes (given by considering $\widehat\Omega$ to be a nonempty bounded subset of an apartment of the building, or of the enlarged building,  of $H(K)$) by extension of scalars $\cO\hookrightarrow \hcO$.   

\vskip.5mm

Let us assume now that Bruhat-Tits theory is available for $H$ over $\hK$. Then Bruhat-Tits theory is also available for $\sD(H)$ over $\hK$\,[P2,\,1.11]. The action of $H(\hK)$ on its building $\cB(\sD(H)/\hK)$ restricts to an action of $H(K)$ by isometries. Let $T$ be a maximal $K$-split torus of $G$ and $A$ be the apartment of $\cB(\sD(H)/\hK)$ corresponding to $T_{\hK}$. We consider the polysimplicial complex $\cB(\sD(H)/\hK)$, with apartments $h\cdot A$, $h\in H(K)$, as the building of $H(K)$ and denote it by $\cB(\sD(H)/K)$. The enlarged Bruhat-Tits building of $H(K)$ is by definition the product $V(\fZ)\times \cB(\sD(H)/K)$.  
\vskip.5mm

Let $\widehat\Omega$ be a nonempty bounded subset of the apartment ${\widehat{A}}= \hh\cdot A$, $\hh\in H(\hK)$, of the building $\cB(\sD(H)/\hK)$, or of the enlarged building $\cB(H/\hK)$. As $H(K)$ is dense in $H(\hK)$, the intersection $\sH_{\widehat\Omega}(\hcO)\hh\cap H(K)$ is nonempty. For any $h$ in this intersection, $\widehat\Omega$ is contained in the apartment $h\cdot A$ of $\cB(\sD(H)/K)$ (or of $\cB(H/K)$). This implies, in particular, that any two facets  lie on an apartment of $\cB(\sD(H)/K)$.  
We now note that the $\hcO$-group schemes $\sH_{\widehat\Omega}$ and $\sH^{\circ}_{\widehat\Omega}$ admit unique descents to smooth affine $\cO$-group schemes with generic fiber $H$, [BLR, Prop.\,D.4(b) in \S6.1]; the affine rings of these descents are $K[H]\cap \hcO[\sH_{\widehat\Omega}]$ and  $K[H]\cap \hcO[\sH^{\circ}_{\widehat\Omega}]$ respectively. \hfill$\Box$

%
\vskip3mm

\ni{\bf 2. Fixed points in $\cB(H/K)$ under a finite automorphism group $\Theta$ of $H$}  
\vskip3mm

We will henceforth assume  that Bruhat-Tits theory is available for $H$ over $K$. 

\vskip2mm

\ni{\bf 2.1.} Let $G$ be a smooth affine $K$-group and $\sG$ be a smooth affine $\cO$-group scheme with generic fiber $G$. According to  [BrT2,\,1.7.1-1.7.2]  $\sG$  is ``\'etoff\'e'' and hence by (ET) of [BrT2,\,1.7.1] its affine ring has the following description:$$\cO[\sG] = \{ f\in K[G]\,\vert\,  f(\sG(\cO))\subset \cO\}.$$ 
\vskip.5mm

Let $\Omega$ be a nonempty bounded subset of an apartment of $\cB(\sD(H)/K)$, or of $\cB(H/K)$.   As the $\cO$-group scheme $\sH_{\Omega}$ is smooth and affine and its generic fiber is $H$, the affine ring of $\sH_{\Omega}$ has thus the  following description:$$\cO[\sH_{\Omega}] = \{ f\in K[H]\,\vert\,  f(H(K)^{\Omega})\subset \cO\}.$$ 
\vskip1mm
 
\ni{\bf Proposition 2.2.} {\em Let $\Omega$ be a nonempty bounded subset of an apartment of $\cB(H/K)$, or of $\cB(\sD(H)/K)$. Let $\sH_{\Omega}$ and $\sH^{\circ}_{\Omega}$ be as above. Let $G$ be a smooth connected $K$-subgroup of $H$ and $\sG$ be a 
smooth affine $\cO$-group scheme with generic fiber $G$ and connected special fiber. Assume that a subgroup $\cG$ of \,$\sG(\cO)$ of finite index fixes $\Omega$ pointwise $($so, $\cG\subset H(K)^{\Omega})$. Then there is a $\cO$-group scheme homomorphism $\varphi:\, \sG\rightarrow \sH^{\circ}_{\Omega}$ that is the natural inclusion $G\hookrightarrow H$ on the generic fibers. So the subgroup  $\sG(\cO)$ of $G(K)$ is contained in $\sH^{\circ}_{\Omega}(\cO)$ and hence it fixes $\Omega$ pointwise. If $F$ is a facet of $\cB(\sD(H)/K)$ that meets $\Omega$, then $\sG(\cO)$  fixes $F$ pointwise. 
\vskip1mm

Let $S$ be a $K$-split torus of $H$ and $\sS$ the split $\cO$-torus with generic fiber $S$. If  a subgroup of the maximal bounded subgroup $\sS(\cO)$ of $S(K)$ of finite index fixes $\Omega$ pointwise, then there is a maximal $K$-split torus $T$ of $H$ containing $S$ such that $\Omega$ is contained in the apartment corresponding to $T$. }    
\vskip2mm

\ni{\it Proof.} Since the fibers of the smooth affine group scheme $\sG$ are connected and the residue field $\kappa$ is separably closed, the subgroup $\cG$ is Zariski-dense in $G$, and its image in $\sG(\kappa)$ is Zariski-dense in the spacial fiber of $\sG$. Using this observation, we easily see that the affine ring $\cO[\sG]\,(\subset K[G])$ of $\sG$ has the following description (cf.\, [BrT2,\,1.7.2]):  $$\cO[\sG] = \{ f\in K[G]\,\vert\,  f(\cG)\subset \cO\}.$$  This description of $\cO[\sG]$ implies at once that the inclusion $\cG\hookrightarrow H(K)^{\Omega}$ induces a $\cO$-group scheme homomorphism  $\varphi:\,\sG\rightarrow \sH_{\Omega}$ that is the natural inclusion $G\hookrightarrow H$ on the generic fibers\,(cf.\,[BrT2, Prop.\,1.7.6]). Since $\sG$ has connected fibers, the homomorphism 
$\varphi$ factors through $\sH^{\circ}_{\Omega}$. 
\vskip.5mm 

Any facet $F$ of $\cB(\sD(H)/K)$ that meets $\Omega$ is stable under $\sG(\cO)\,(\subset H(K))$, so a subgroup of $\sG(\cO)$ of finite index fixes it pointwise. Now applying the result of the preceding paragraph, for $F$ in place of $\Omega$, we see that there is a $\cO$-group scheme homomorphism $\sG\rightarrow \sH^{\circ}_F$ that is the natural inclusion $G\hookrightarrow H$ on the generic fibers and hence $\sG(\cO)$ fixes $F$ pointwise. 
\vskip.5mm

Now we will prove the last assertion of the proposition. It follows from what we have shown above  that there is a $\cO$-group scheme homomorphism $\iota: \sS\rightarrow \sH^{\circ}_{\Omega}$ that is the natural inclusion $S\hookrightarrow H$ on the generic fibers ($\iota$ is actually a closed immersion, see [SGA3$_{\rm{II}}$, Exp.\,IX, Thm.\,6.8 and Cor.\,2.5] and also the proof of Lemma 4.1 in [PY2]). Applying [P2, Prop.\,2.1(i)] to the centralizer of $\iota(\sS)$ (in $\sH^{\circ}_{\Omega}$) in place of $\sG$, and $\cO$ in place of $\fo$, we see that there is a closed $\cO$-torus $\sT$ of $\sH^{\circ}_{\Omega}$ that commutes with  $\iota(\sS)$ and whose generic fiber $T$ is a maximal $K$-split torus of $H$. The torus $T$ clearly contains $S$, and [P2,\,Prop.\,2.2(ii)] implies that $\Omega$ is contained in the apartment corresponding to $T$.        \hfill$\Box$
\vskip2mm

The following is a simple consequence of the preceding proposition.
\vskip1mm

\ni{\bf Corollary 2.3.} {\em Let $G$, $S$, $\sG$, and $\sS$ be as in the preceding proposition. Then the set of points of $\cB(\sD(H)/K)$ that are fixed under $\sG(\cO)$  is the union of facets pointwise fixed under $\sG(\cO)$. The set of points of the enlarged building $\cB(H/K)$ that are fixed under a finite-index subgroup $\cS$ of the maximal bounded subgroup $S(K)_b\,(=\sS(\cO))$ of $S(K)$  is  the enlarged Bruhat-Tits building  $\cB(Z_H(S)/K)$ of the centralizer $Z_H(S)(K)$ of $S$ in $H(K)$.}
\vskip2mm
{
\vskip2mm
\ni{\bf 2.4.} Let $\Theta$ be a finite group of automorphisms of  the reductive  $K$-group $H$. There is a natural action of $\Theta$ on $\cB(H/K)$, as well as on $\cB(\sD(H)/K)$, by  isometries such that for all $h$  in $H(K)$, $x$ in  $\cB(H/K)$ or $\cB(\sD(H)/K)$,  and $\theta$ in  $\Theta$, we have $\theta(h\cdot x)= \theta(h)\cdot \theta(x)$.  The action of $\Theta$ on $\cB(\sD(H)/K)$ is polysimplicial. 
\vskip1mm

In the following we will assume that the characteristic $p$  of the residue field $\kappa$ does {\it not} divide the order of $\Theta$. Then $G := (H^{\Theta})^{\circ}$ is a reductive group, see [Ri, Prop.\,10.1.5] or [PY1, Thm.\,2.1]. We will prove that Bruhat-Tits theory is available for $G$ over $K$ and the enlarged Bruhat-Tits building of $G(K)$, as a metric space, can be identified with the closed subspace $\cB(H/K)^{\Theta}$ of points  of $\cB(H/K)$ fixed under $\Theta$. 

\vskip1mm

For a subset $X$ of a set given with an action of $\Theta$, $X^{\Theta}$ will denote the subset of points of $X$ that are fixed under $\Theta$. We will denote $\cB(\sD(H)/K)^{\Theta}$ by $\sB$ in the sequel. 
\vskip1mm

If a facet of $\cB(\sD(H)/K)$ is stable under the action of $\Theta$, then its barycenter is fixed under $\Theta$. Conversely, if a facet $F$ of $\cB(\sD(H)/K)$ contains a point $x$ fixed under $\Theta$, then being the unique facet containing $x$, $F$ is stable under the action of $\Theta$. 

\vskip2mm

 \ni{\bf 2.5}. (i) {\bf Bruhat-Tits group schemes $\sG^{\circ}_{\Omega}$.}  Let $\Omega$ be a  nonempty $\Theta$-stable  bounded subset of an apartment of $\cB(\sD(H)/K)$, or of $\cB(H/K)$.  Then $\sH_{\Omega}(\cO)\,(= H(K)^{\Omega})$ is stable under  the action of $\Theta$ on $H(K)$, so the affine ring $\cO[\sH_{\Omega}]$ is stable under the action of $\Theta$  on $K[H]$. This implies that $\Theta$  acts on the group scheme $\sH_{\Omega}$ by $\cO$-group scheme automorphisms. The neutral component  $\sH^{\circ}_{\Omega}$ of $\sH_{\Omega}$  is of course stable under this action.  Define the functor $\sH_{\Omega}^{\Theta}$ of $\Theta$-fixed points that associates to a commutative $\cO$-algebra $C$  the subgroup $\sH_{\Omega}(C)^{\Theta}$ of 
$\sH_{\Omega}(C)$ consisting of elements fixed under $\Theta$. The functor $\sH_{\Omega}^{\Theta}$ is represented by a closed smooth $\cO$-subgroup scheme of $\sH_{\Omega}$ (see Propositions 3.1 and 3.4 of [E], or Proposition A.8.10 of [CGP]); we will denote this closed smooth $\cO$-subgroup scheme also by $\sH_{\Omega}^{\Theta}$.   Its generic fiber is $H^{\Theta}$, and so the identity component of the generic fiber is $G$.  The {\it neutral component} ${(\sH_{\Omega}^{\Theta})}^{\circ}$ of $\sH_{\Omega}^{\Theta}$ is by definition the union of the identity components of its generic and special fibers; it is an open (so smooth) affine $\cO$-subgroup scheme\,[PY2, \S3.5] with generic fiber $G$. The index of the subgroup ${(\sH_{\Omega}^{\Theta})}^{\circ}(\cO)$ in $\sH_{\Omega}^{\Theta}(\cO)$ is known to be finite\,[EGA\,IV$_3$, Cor.\,15.6.5]. It is obvious that ${(\sH_{\Omega}^{\Theta})}^{\circ}= ((\sH^{\circ}_{\Omega})^{\Theta})^{\circ}$. We will denote $(\sH_{\Omega}^{\Theta})^{\circ}$ by $\sG^{\circ}_{\Omega}$ in the sequel and call it {\it the Bruhat-Tits $\cO$-group scheme associated to $G$ and $\Omega$}.  As $\sG^{\circ}_{\Omega}(\cO)\subset\sH_{\Omega}(\cO)$, $\sG^{\circ}_{\Omega}(\cO)$ fixes $\Omega$ pointwise. 
\vskip1mm

(ii) {\bf Parahoric subgroups of $G(K)$.} For $x\in \sB$, $\sG^{\circ}_{\{x\}}$ and $P_x := \sG^{\circ}_{\{x\}}(\cO)$ will  respectively be called the {\it Bruhat-Tits parahoric $\cO$-group scheme} and the {\it parahoric subgroup} of $G(K)$ associated with the point $x$.  The parahoric subgroups of $G(K)$ are by definition the subgroups  $P_x$ for $x\in \sB$. For a given parahoric subgroup $P_x$, the associated Bruhat-Tits parahoric $\cO$-group scheme is $\sG^{\circ}_{\{x\}}$. 
\vskip1mm

(iii) Let $P$ be a parahoric subgroup of $G(K)$, $\sG^{\circ}$ the associated Bruhat-Tits parahoric $\cO$-group scheme, $\osG^{\circ}$ the special fiber of $\sG^{\circ}$,  and $\cP$ be a subgroup of $P$ of finite index. Then the image of $\cP$ in $\osG^{\circ}(\kappa)$ is Zariski-dense in the connected group $\osG^{\circ}$, so  the  affine ring of $\sG^{\circ}$ is: $$\cO[\sG^{\circ}]=\{f\in K[G]\,|\, f(\cP)\subset \cO\}.$$ Thus the subgroup $\cP$ ``determines'' the group scheme $\sG^{\circ}$, and hence $P$ is the unique parahoric subgroup of $G(K)$ containing $\cP$ as a subgroup of finite index.
\vskip1mm

(iv) Let $P$ and $\sG^{\circ}$ be as in the preceding paragraph. Let $\Omega$ be a nonempty $\Theta$-stable bounded subset of an apartment of $\cB(\sD(H)/K)$, or of $\cB(H/K)$, and  $\sG^{\circ}_{\Omega}$ be as in (i). We assume that $\Omega$ is fixed pointwise by $P$. Then the inclusion of $P$ in $H(K)^{\Omega} \,(=\sH_{\Omega}(\cO))$ gives a $\cO$-group scheme homomorphism $\sG^{\circ}\rightarrow \sH^{\circ}_{\Omega}$ (Proposition 2.2). This homomorphism obviously  factors through $\sG^{\circ}_{\Omega}$ to give a $\cO$-group scheme homomorphism $\sG^{\circ}\rightarrow \sG^{\circ}_{\Omega}$ that is the identity on the generic fiber $G$.   
\vskip.5mm

Suppose $x,\,y$ are two points of $\sB$  that are fixed by $P$, and $[xy]$ is the geodesic joining $x$ and $y$. Then $P$ fixes every point $z$ of $[xy]$. Let $\sG^{\circ}_{[xy]}$ be as in (i) for $\Omega=[xy]$. There are $\cO$-group scheme homomorphisms $\sG^{\circ}\rightarrow \sG^{\circ}_{[xy]}$ and $\sG^{\circ}\rightarrow \sG^{\circ}_{\{z\}}$ that are the identity on the generic fiber $G$.  
\vskip2mm 

We will show that $\sB$ is an affine building with apartments described in 2.16. We begin with the following proposition which has been suggested by Proposition 1.1 of  [PY1]. The proof given here is different from that in [PY1]. 
\vskip2mm

\ni{\bf Proposition 2.6.} {\em Let $\Omega$ be a nonempty $\Theta$-stable bounded subset of an apartment of $\cB(\sD(H)/K)$. Let $\sH :=\sH^{\circ}_{\Omega}$ be the Bruhat-Tits smooth affine $\cO$-group scheme with generic fiber $H$, and connected special fiber $\osH$, associated to  $\Omega$.   Let $\orH : = \osH/{\sR_{u,\kappa}}(\osH)$ be the maximal pseudo-reductive quotient of $\osH$.  Then there exist closed $\cO$-tori $\sS\subset \sT$ in $\sH$ such that
\vskip1mm

{\rm(i)} the generic fiber $T$ of $\sT$ is a maximal $K$-split torus of  $H$ and $\Omega$ is contained in the apartment corresponding to $T$;
\vskip1mm

{\rm(ii)} $\sS$, and therefore its generic fiber $S$, is stable under $\Theta$ and the special fiber of $\sS$ maps onto the  central torus  of $\orH$. }
\vskip2mm

%
 
 We first prove the following lemma in which $\kappa$ is any field of characteristic $p\geqslant 0$.
 \vskip2mm
 
 \ni{\bf Lemma 2.7.} {\em Let $\cH$ be a smooth connected affine algebraic $\kappa$-group given with an action by a finite group $\Theta$ and $\cU$ be a $\Theta$-stable smooth connected  unipotent normal $\kappa$-subgroup of $\cH$. We assume that $p$ does not divide  the order of $\Theta$. Let $\ocS$ be a $\Theta$-stable $\kappa$-torus of $\ocH:= \cH/\cU$. Then there exists a $\Theta$-stable $\kappa$-torus $\cS$ in $\cH$ that maps isomorphically onto $\ocS$. In particular, there exists a $\Theta$-stable $\kappa$-torus in $\cH$ that maps isomorphically onto the central torus of $\ocH$.}
\vskip2mm 

\ni{\it Proof.} Let $\cT$ be a $\kappa$-torus of $\cH$ that maps isomorphically onto $\ocS\,(\subset \ocH)$. Considering the $\Theta$-stable solvable subgroup $\cT\,\cU$; using conjugacy under $\cU(\kappa)$ of maximal $\kappa$-tori of this solvable group [Bo,\:Thm.\,19.2], we see that for $\theta\in \Theta$,  $\theta(\cT) = u(\theta)^{-1}\cT u(\theta)$ for some $u(\theta) \in \cU(\kappa)$.   Let $\cU(\kappa)=: \cU_0\supset \cU_1\supset \cU_2\cdots \supset \cU_n =\{1\}$ be the descending central series of the nilpotent group $\cU(\kappa)$. Each subgroup $\cU_i$ is $\Theta$-stable and $\cU_i/\cU_{i+1}$ is a commutative $p$-group if $p\ne 0$, and a $\bQ$-vector space if $p =0$. Now let $i\leqslant n$, be the largest integer such that there exists a $\kappa$-torus $\cS$  in $\cT\,\cU$ that maps onto $\ocS$, and for every $\theta\in \Theta$, there is a $u(\theta)\in \cU_i$ such that $\theta(\cS) = u(\theta)^{-1}\cS u(\theta)$. Let $\cN_i$ be the normalizer of $\cS$ in $\cU_i$. Then, for $\theta\in \Theta$, $\theta (\cN_i) = u(\theta)^{-1}\cN_i u(\theta)$ and hence as $\cU_i/\cU_{i+1}$ is commutative, we see that $\theta(\cN_i \cU_{i+1}) = \cN_i \cU_{i+1}$, i.e., $\cN_i \cU_{i+1}$ is $\Theta$-stable.  It is easy to see that $\theta\mapsto u(\theta)$ mod\,$(\cN_i \cU_{i+1})$ is a $1$-cocycle on $\Theta$ with values in $\cU_i/\cN_i \cU_{i+1}$. But ${\rm{H}}^1(\Theta,\, \cU_i/\cN_i \cU_{i+1})$ is trivial since the finite group $\Theta$ is of order prime to $p$ if $p\ne 0$, and $\cU_i/\cN_i \cU_{i+1}$ is divisible if $p =0$.   So there exits a $u\in \cU_i$ such that for all $\theta\in\Theta$, $u^{-1}u(\theta)\theta(u)$ lies in $\cN_i \cU_{i+1}$.  Now let  $\cS' =u^{-1}\cS u$. Then the normalizer of $\cS'$ in $\cU_i$ is $u^{-1} \cN_i u$ and again as $\cU_i/\cU_{i+1}$ is commutative, $u^{-1}\cN_i u\cdot \cU_{i+1} = \cN_i \cU_{i+1}$. For $\theta \in \Theta$, we choose  $u'(\theta)\in \cU_{i+1}$ such that $u^{-1}u(\theta)\theta(u)\in u^{-1}\cN_i u\cdot u'(\theta)$.  Then $\theta(\cS') = u'(\theta)^{-1}\cS' u'(\theta)$ for all $\theta \in \Theta$. This contradicts the maximality of $i$ unless $i =n$. \hfill$\Box$ 
\vskip2mm

 We will now prove Proposition 2.6. There is a natural action of $\Theta$ on $\sH$ by $\cO$-group scheme automorphisms\,(2.5(i)). This action induces an action of $\Theta$ on $\osH$, as well as on $\osH^{\rm{pred}}$, by $\kappa$-group automorphisms.  According to the preceding lemma, there exists a $\Theta$-stable $\kappa$-torus $\osS$ in $\osH$  that maps isomorphically onto the central torus of $\osH^{\rm{pred}}$. We will denote the natural inclusion $\osS\hookrightarrow \osH$ by $\overline\iota$. The character group ${\rm{X}}^{\ast}(\osS)$ of $\osS$ is a free abelian group on which $\Theta$ acts by automorphisms. Let $\sS$ be the split $\cO$-torus with affine ring the group ring $\cO[{\rm{X}}^{\ast}(\osS)]$. The action of $\Theta$ on  ${\rm{X}}^{\ast}(\osS)$ induces an action (of $\Theta$) on $\sS$ by $\cO$-group scheme automorphisms. According to [SGA3$_{\rm{II}}$, Exp.\,XI,\,4.2] the homomorphism functor ${\underline{\rm{Hom}}}_{{\rm {Spec}}(\cO){{\textendash}\rm{gr}}}(\sS,\sH)$   is representable by a smooth $\cO$-scheme $\sX$. There is a natural action of $\Theta$ on $\sX$ induced by its action on  ${\underline{\rm{Hom}}}_{{\rm {Spec}}(\cO){{\textendash}\rm{gr}}}(\sS,\sH)$ described as follows. $$ {\rm {For}}\: f\in {\underline{\rm{Hom}}}_{{\rm {Spec}}(\cO){{\textendash}\rm{gr}}}(\sS,\sH),\:  \theta\in \Theta, \: {\rm{and}}\: s\in \sS(C),\: (\theta\cdot f)(s) = \theta(f(\theta^{-1}(s))),$$ for any commutative $\cO$-algebra $C$.  By Proposition A.8.10 of [CGP] the subscheme $\sX^{\Theta}$ of $\Theta$-fixed points in $\sX$ is a closed smooth $\cO$-subscheme of $\sX$.   The inclusion $\overline{\iota}:\, \osS\hookrightarrow \osH$ clearly lies in $\sX^{\Theta}(\kappa)$. Now since $\cO$ is Henselian, the natural map $\sX^{\Theta}(\cO)\rightarrow \sX^{\Theta}(\kappa)$ is surjective\,[EGA\,IV$_4$,\,18.5.17] and hence there is a $\cO$-group scheme homomorphism $\iota:\, \sS\rightarrow\sH$ that is $\Theta$-equivariant and lies over $\overline{\iota}$. As $\overline{\iota}$ is a closed immersion, using [SGA3$_{\rm{II}}$, Exp.\,IX, 2.5 and 6.6] we see that $\iota$ is a closed immersion. We identify $\sS$ with a $\Theta$-stable closed $\cO$-torus of $\sH$ in terms of $\iota$.  
 \vskip1mm
 
 We now consider the centralizer $\mathscr{C}$ of $\sS$ in $\sH$ which by [CGP, Prop.\,A.8.10(2)] is a smooth affine $\cO$-subgroup scheme of $\sH$.  Its special fiber clearly contains a maximal torus $\osT$ of $\osH$. Hence, using Proposition 2.1(i) of [P2] we see that $\mathscr{C}$ contains a closed $\cO$-torus $\sT$ whose special fiber is $\osT$. This $\sT$ is $\cO$-split since its special fiber $\osT$ is split over the separably closed field $\kappa$; $\sT$ clearly contains $\sS$, and its generic fiber $T$ is a maximal $K$-split torus of $G$ since $\osT$ is a maximal $\kappa$-torus of $\osH$. Proposition 2.2(ii) of [P2] implies that the apartment of $\cB(\sD(H)/K)$ corresponding to $T$ contains $\Omega$. \hfill$\Box$   

\vskip2mm
\ni{\bf 2.8. A reduction.} Let $Z$ be the central torus of $H$, then $H = Z\cdot\sD(H)$. Hence, $G=(H^{\Theta})^{\circ} = (Z^{\Theta})^{\circ}\cdot(\sD(H)^{\Theta})^{\circ}$.  
As the enlarged Bruhat-Tits building $\cB(H/K)$ of $H(K)$ is by definition $V(\fZ)\times \cB(\sD(H)/K)$,  $\cB(H/K)^{\Theta} = V(\fZ)^{\Theta}\times \cB(\sD(H)/K)^{\Theta}$.  So to prove that $\cB(H/K)^{\Theta}$, as a metric space, is the enlarged Bruhat-Tits building of $G(K)$, it would suffice to show that as a metric space, $\cB(\sD(H)/K)^{\Theta}$ is the enlarged Bruhat-Tits building of $(\sD(H)^{\Theta})^{\circ}$.   In view of this, we may (and do) replace $H$ by $\sD(H)$ and assume in the sequel that $H$ {\it is semi-simple}.   
\vskip2mm

\ni{\bf 2.9.}  We introduce the following partial order ``$\prec$'' on the set of nonempty subsets of $\cB(H/K)$: Given two nonempty subsets $\Omega$ and $\Omega'$, $\Omega'\prec \Omega$ if the closure $\overline{\Omega}$ of $\Omega$ contains $\Omega'$. If $F$ and $F'$ are facets of $\cB(H/K)$, with $F'\prec F$, or equivalently, $\sH^{\circ}_{F}(\cO)\subset \sH^{\circ}_{F'}(\cO)$, we say that  $F'$ is a {\it face} of $F$.  In a collection $\cC$ of facets, thus a facet is {\it maximal} if it is not a proper face of any facet belonging to  $\cC$, and a facet is {\it minimal} if no proper face of it belongs to $\cC$.  
\vskip1mm

Now let $X$ be a convex subset of $\cB(H/K)$ and $\cC$ be the set of facets of $\cB(H/K)$, or facets lying in a given apartment $A$, that meet $X$. Then the following assertions are easy to prove\,(see Proposition 9.2.5 of [BrT1]): (1) All maximal facets in $\cC$ are of equal dimension and a facet $F\in \cC$ is maximal if and only if $\dim(F\cap X)$ is maximal.   (2) Let $F$ be a facet lying in an apartment $A$. Assume that $F$ is  maximal among the facets of $A$ that meet $X$, and let $A_F$ be the affine subspace of $A$ spanned by $F$. Then every facet of $A$ that meets $X$ is contained in $A_F$ and $A\cap X$ is contained in the affine subspace of $A$ spanned by $F\cap X$. 
 \vskip1.5mm

The subset $\sB \,(=\cB(H/K)^{\Theta})$ of $\cB(H/K)$ is closed and convex. Hence the assertions of the preceding paragraph hold for $\sB$ in place of $X$. 
\vskip2mm

\ni{\bf 2.10.}  Let $x, \,y\in \sB \,(= \cB(H/K)^{\Theta})$. Let $F$ be a facet of $\cB(H/K)$ which contains $x$ in its closure and  
is maximal among the facets that meet $\sB$, and let $\Omega = F\cup \{y\}$.  Let  $S\subset T$ be a pair of $K$-split tori with properties (i) and (ii) of Proposition 2.6.  Let $S_G$ and $T_G$ be the maximal subtori of $S$ and $T$ respectively contained in $G$.  Let $A$ be the apartment of $\cB(H/K)$ corresponding to $T$. Then $A$ contains $y$ and the closure of $F$, and so it also contains $x$. Moreover, $A$ is an affine space under $V(T)$, the affine subspace $V(S)+x$ of $A$ contains $F$ and is spanned by it. The affine subspaces  $V(S_G)+x \subset V(T_G)+x$ of $A$ are  clearly contained in $\sB$. As $V(S)^{\Theta} = V(S_G)$ and  $F\subset V(S)+x$, we see that  $F^{\Theta}= F\cap \sB$ is contained in $V(S_G)+x$.   But since the facet $F$ is maximal among the facets that meet  $\sB$, $A^{\Theta}\,(=A\cap \sB)$ is contained in the affine subspace  of $A$ spanned by 
$F^{\Theta}$. Therefore, $A^{\Theta}=V(S_G)+x$. This implies that $V(S_G)+x = V(T_G)+x$ and hence $S_G =T_G$. We will now show that $S_G$ is a maximal $K$-split torus of $G$.    
 \vskip1mm
 
 Let $S'$ be a maximal $K$-split torus of $G$ containing $S_G$.  The centralizer  $Z_H(S')$ of $S'$ in $H$ is stable under $\Theta$. The enlarged Bruhat-Tits building $\cB(Z_H(S')/K)$ of $Z_H(S')(K)$ is identified with the union of apartments of $\cB(H/K)$ corresponding to maximal $K$-split tori of $H$ that contain $S'$, cf.\,[P2, 3.11]. Let $z$ be a point of $\cB(Z_H(S')/K)^{\Theta}$ and $T'\, (\supset S')$ be a maximal $K$-split torus of $H$ such that the corresponding apartment $A'$ of $\cB(Z_H(S')/K)$ contains $z$.   Then $A'=V(T')+z$ and hence ${A'}^{\Theta}=A'\cap \sB=  V(T')^{\Theta}+z =V(S')+z$ is an affine subspace of $A'$ of dimension  $\dim(S')$. Let $F'$ be a facet of  $A'$ that contains the point $z$ in its closure and is maximal among the facets of $A'$ meeting $\sB$. Then ${A'}^{\Theta}$ is contained in the affine subspace of $A'$ spanned by ${F'}^{\Theta}$, so $\dim({F'}^{\Theta}) = \dim(S')\geqslant \dim(S_G)$. But $\dim(F^{\Theta}) = \dim(S_G)\geqslant \dim({F'}^{\Theta})$. This implies that $\dim(S_G)=\dim(S')$ and hence $S' = S_G$. Thus we have shown that  $S_G$ is a maximal $K$-split torus of $G$.  
 
 \vskip1mm
 
 Now let $x$ be an arbitrary point of $\sB$ and take $y =x$. Let $S_G$ be as above and $C$ be the maximal $K$-split  central torus of $G$. Then $C$ is contained in $S_G$ and hence  $Z_H(S_G)$ is contained in the centralizer $Z_H(C)$ of $C$ in $H$. So $\cB(Z_H(S_G)/K)$ is contained in $\cB(Z_H(C)/K)$.  In particular,  $x\in \cB(Z_H(C)/K)$, which implies that   $\sB\subset \cB(Z_H(C)/K)$. 
 \vskip2mm
 
 Thus we have established the following proposition:
 \vskip2mm
 
 \ni{\bf Proposition 2.11.} (i) {\em  $\sB\subset \cB(Z_H(C)/K)$.}  
 \vskip.5mm
 
         (ii) {\em Given points $x,\,y\in \sB$, there exists a maximal $K$-split torus $S_G$ of $G$, and a maximal $K$-split torus $T$ of $H$ containing $S_G$, and hence contained in $Z_H(S_G)$, such that the apartment $A\,(\subset \cB(Z_H(S_G)/K))$ corresponding to $T$ contains $x$ and $y$. Moreover, $A^{\Theta}=A\cap \sB$ is the affine subspace $V(S_G)+x$ of $A$ of dimension $\dim(S_G)$.}
 \vskip2mm
 
 We will now derive the following proposition which will give us apartments in the enlarged Bruhat-Tits building of $G(K)$. In the sequel, we will use  $S$, instead of $S_G$,  to denote a maximal $K$-split torus of $G$. As $M:=Z_H(S)$ is stable under $\Theta$, the enlarged Bruhat-Tits building  $\cB(M/K)$ of $M(K)$ contains a $\Theta$-fixed point. 
 \vskip2mm
 
 \ni{\bf Proposition 2.12.} {\em Let $S$ be a maximal $K$-split torus of $G$ and let $T$ be a maximal $K$-split torus of $H$ containing $S$ such that the apartment $A$ of $\cB(H/K)$ corresponding to $T$ contains a $\Theta$-fixed point $x$. Then $\cB(Z_H(S)/K)^{\Theta} = V(S)+x = A^{\Theta}$. So $\cB(Z_H(S)/K)^{\Theta}$ is an affine space under the $\bR$-vector space $V(S)$.}  
\vskip2mm

\ni{\it Proof.}  Let $Z$ be the maximal $K$-split central torus and $Z_H(S)'$  be the derived subgroup of $Z_H(S)$. Then $Z$, $Z_H(S)$ and $Z_H(S)'$ are stable under $\Theta$; $G':={({Z_H(S)'}^{\Theta})}^{\circ}$ is anisotropic over $K$ since $S$ is a maximal $K$-split torus of $(Z_H(S)^{\Theta})^{\circ}\,(\subset G)$.  Now applying Proposition 2.11(ii)  to $Z_H (S)'$ in place of $H$, we see that the Bruhat-Tits building $\cB(Z_H(S)'/K)$ of $Z_H(S)'(K)$ contains only one point fixed under $\Theta$. For if $y, z\in\cB(Z_H(S)'/K)^{\Theta}$, then there is an apartment $A'$ of $\cB(Z_H(S)'/K)$ containing these points and ${A'}^{\Theta}$  is an affine subspace of $A'$  of dimension 0 as $G'$ is anisotropic over $K$. Therefore, $y=z$.     This proves that $\cB(Z_H(S)'/K)^{\Theta}$ consists of a single point. Hence, $\cB(Z_H(S)/K)^{\Theta} = V(Z)^{\Theta}+x = V(S)+x$, and so it is an affine space under $V(S)$. \hfill$\Box$

\vskip2mm

\ni{\bf 2.13. Another reduction.} Let $C$ be the maximal $K$-split  central torus of $G$. Then $\cB(H/K)^{\Theta}=\sB\subset \cB(Z_H(C)/K)$ (Proposition 2.11(i)).  Let $H'$ be the derived subgroup of $Z_H(C)$; $H'$ is a connected semi-simple subgroup of $H$ stable under the group $\Theta$ of automorphisms of $H$; ${({H'}^{\Theta})}^{\circ}\,(\subset G)$ contains the derived subgroup of $G$ and its central torus is $K$-anisotropic.  To prove that $\cB(H/K)^{\Theta}$ is the enlarged Bruhat-Tits building of $G(K)$, {\em we} {may (and} {\em will}) {\em assume, after replacing $H$ with $H'$, that $H$ is semi-simple and the central torus of $G$ is $K$-anisotropic {\rm(cf.\,[P2, 3.11,\,1.11])} and show that $\cB(H/K)^{\Theta}$ is an affine building}. 

\vskip2mm

\ni{\bf 2.14.} Let $S$ be a maximal $K$-split torus of $G$. Let $N :=N_G(S)$ and $Z := Z_G(S)$ respectively  be the normalizer and the centralizer of $S$ in $G$. As $N$ (in fact, the normalizer $N_H(S)$ of $S$ in $H$) normalizes the centralizer $Z_H(S)$ of $S$ in $H$, there is a natural action of  $N(K)$ on $\cB(Z_H(S)/K)$ and $N(K)$ 
stabilizes $\cB(Z_H(S)/K)^{\Theta}$ under this action. For $n\in N(K)$, the action of $n$ carries an apartment $A$ of $\cB(Z_H(S)/K)$  to the apartment $n\cdot A$ by an affine transformation. 
\vskip1mm

 Now let $T$ be a maximal $K$-split torus of $Z_H(S)$ such that the corresponding apartment $A:=A_T$ of $\cB(Z_H(S)/K)$ contains a $\Theta$-fixed point $x$.  According to the previous proposition, $\cB(Z_H(S)/K)^{\Theta} =  V(S)+x= A^{\Theta}$. So we can view $\cB(Z_H(S)/K)^{\Theta}$ as an affine space under $V(S)$. We will now show, using the proof of the lemma in 1.6 of [PY1], that $\cB(Z_H(S)/K)^{\Theta}$ has the properties required of an apartment corresponding to the maximal $K$-split torus $S$ in the Bruhat-Tits building of $G(K)$ if such a building exists.  We need to check the following three conditions.  
 \vskip1.5mm
 
\ni{A1:} {\em The action of $N(K)$ on $\cB(Z_H(S)/K)^{\Theta}= A^{\Theta}$ is by affine transformations and the maximal bounded subgroup $Z(K)_b$ of $Z(K)$ acts trivially.}
\vskip1mm

Let ${\rm Aff}(A^{\Theta})$ be the group of affine automorphisms of $A^{\Theta}$ and $\varphi: N(K)\rightarrow {\rm Aff}(A^{\Theta})$ be the action map.
\vskip1mm
 
\ni{A2:} {\em The group $Z(K)$ acts by translations, and the action is characterized by the following 
formula: for $z \in Z(K)$,
$$\chi(\varphi(z)) = -\omega(\chi(z)) \ {\rm for \ all}\  \chi \in {{\rm X}_K^{\ast}}(Z)\,(\hookrightarrow {{\rm X}_K^{\ast}}(S)),$$
here we regard the translation $\varphi(z)$ as an element of  $V(S)$}.
\vskip1mm

\ni{A3:} {\em For $g \in {\rm{Aff}}(A^{\Theta})$, denote by $dg \in {\rm{GL}}(V(S))$ the derivative of $g$. Then the map 
$N(K)\rightarrow  {\rm{GL}}(V(S))$, $n\mapsto d\varphi(n)$, is induced from the action of $N(K)$ on ${{\rm X}_{\ast}}(S)$ $($i.e., it is the Weyl group action$)$.}
\vskip1mm

Moreover, as the central torus of $G$ is $K$-anisotropic, these three conditions determine the affine structure on $\cB(Z_H(S)/K)^{\Theta}$ uniquely; see [T, 1.2].
\vskip2mm

\ni{\bf Proposition 2.15.} {\em Conditions {\rm {A1, A2}} and {\rm{A3}} hold.}
\vskip2mm

\ni{\it Proof.} The action of $n \in N(K)$ on $\cB(Z_H(S)/K)$ carries the apartment $A = A_T$  via an affine isomorphism $f(n)\,:\, A\rightarrow A_{nTn^{-1}}$ to the apartment $A_{nTn^{-1}}$ corresponding to the torus $nTn^{-1}$ containing $S$. As $(A_{nTn^{-1}})^{\Theta} = \cB(Z_H(S)/K)^{\Theta} = A^{\Theta}$, we see that $f(n)$ keeps $A^{\Theta}$ stable and so $\varphi(n):= f(n)\vert_{A^{\Theta}}$ is an affine automorphism of $A^{\Theta}$.  
\vskip1mm

The derivative $df(n): V(T) \rightarrow V(nTn^{-1})$  is induced from the map $${\rm {Hom}}_K({\rm {GL}}_1, T)={{\rm X}_{\ast}}(T)\rightarrow  {{\rm X}_{\ast}}(nTn^{-1})={\rm{Hom}}_K({\rm{GL}}_1, nTn^{-1}),$$ $\lambda\mapsto {\rm{Int}}\,n\cdot\lambda$, where ${\rm{Int}}\,n$ is the inner automorphism of $H$ determined by $n\in N(K)\subset H(K)$. So, the restriction $d\varphi(n):V(S)\rightarrow V(S)$ is induced from the homomorphism ${\rm X}_{\ast}(S)\rightarrow {\rm X}_{\ast}(S)$, $\lambda\mapsto {\rm{Int}}\,n\cdot \lambda$. This proves {A3}.

\vskip1mm
Condition {A3} implies that $d\varphi$ is trivial on $Z(K)$. Therefore, $Z(K)$ acts by translations. 
The action of the bounded subgroup $Z(K)_b$ on $A^{\Theta}$ admits a fixed point by the fixed point theorem of Bruhat-Tits. Therefore, $Z(K)_b$ acts by the trivial translation. This proves {A1}.
\vskip1mm

Since the image of $S(K)$ in $Z(K)/Z(K)_b\simeq \bZ^{\dim(S)}$ is a subgroup of finite index, to prove the formula in {A2}, it suffices to prove it for $z\in S(K)$. But for $z\in S(K)$, $zTz^{-1} = T$, and $f(z)$ is a translation of the apartment $A$ ($\varphi(z)$ is regarded as an element of $V(T)$) which satisfies\,(see 1.9 of [P2]):  $$\chi(f(z)) =-\omega(\chi(z))\ \ {\rm for \ all} \ \chi\in {\rm X}_K^{\ast}(T).$$ This implies the formula in {A2}, since the restriction map ${\rm X}_K^{\ast}(T)\rightarrow {\rm X}_K^{\ast}(S)$ is surjective and the image of the restriction map ${\rm X}_K^{\ast}(Z)\rightarrow {\rm X}_K^{\ast}(S)$ is of finite index in ${\rm X}_K^{\ast}(S)$. \hfill$\Box$
\vskip3mm

\ni{\bf 2.16. Apartments of $\sB$.} By definition, the {\it apartments} of  $\sB$ are the affine spaces $\cB(Z_H(S)/K)^{\Theta}$ under the $\bR$-vector space $V(S)$ (of dimension $=K$-rank\,$G$) for maximal $K$-split tori $S$ of $G$. For any apartment $A$ of $\cB(Z_H(S)/K)$ that contains a $\Theta$-fixed point, $\cB(Z_H(S)/K)^{\Theta} = A^{\Theta}$\,(Proposition 2.11(ii)). The subgroup $N_G(S)(K)$ of $G(K)$ acts by affine transformations on the apartment $\cB(Z_H(S)/K)^{\Theta}$  and $Z_G(S)(K)$ acts on it by translations (Proposition 2.15).  Conjugacy of maximal $K$-split tori of $G$ under $G(K)$ implies that this group acts transitively on the set of apartments of $\sB$. 
\vskip1mm

Propositions 2.11(ii) and 2.12 imply the following proposition at once:
\vskip2mm

\ni{\bf Proposition 2.17.} {\em Given  any two points of $\sB$,  there is a maximal $K$-split torus $S$ of $G$ such that the corresponding apartment of $\sB$ contains these two points.}
\vskip2mm 

\ni{\bf Proposition 2.18.} {\em Let $\cA$ be an apartment of $\sB$. Then there is a {unique} maximal $K$-split torus $S$ of $G$ 
such that $\cA = \cB(Z_H(S)/K)^{\Theta}$. So the stabilizer of $\cA$ in $G(K)$ is $N_G(S)(K)$.} 
\vskip2mm

\ni{\it Proof.}  We fix a maximal $K$-split torus $S$ of $G$ such that $\cA = \cB(Z_H(S)/K)^{\Theta}$. We will show that $S$ is uniquely determined by $\cA$. For this purpose, we observe that the subgroup $N_G(S)(K)$ of $G(K)$ acts on $\cA$ and the maximal bounded subgroup $Z_G(S)(K)_b$ of $Z_G(S)(K)$ acts trivially (Proposition 2.15). So the subgroup $\mathcal{Z}$ of $G(K)$ consisting of elements that fix $\cA$ pointwise is a bounded subgroup of $G(K)$, normalized by $N_G(S)(K)$,  and it contains $Z_G(S)(K)_b$. Now, using the Bruhat decomposition of $G(K)$ with respect to $S$, we see that every bounded subgroup of $G(K)$ that is normalized by $N_G(S)(K)$ is a normal subgroup of the latter.  Hence the identity component of the Zariski-closure of $\cZ$ is $Z_G(S)$.  As $S$ is the unique maximal $K$-split torus of $G$ contained in $Z_G(S)$, both the assertions follow.\hfill$\Box$
\vskip2mm

\ni {\bf 2.19. The  affine Weyl group of $G$.}  Let $G(K)^+$ denote the (normal) subgroup of $G(K)$ generated by $K$-rational elements of the unipotent radicals of parabolic $K$-subgroups of $G$. Let $S$ be a maximal $K$-split torus of $G$, $N$ and $Z$ respectively be the normalizer and centralizer of $S$ in $G$. Let $N(K)^{+} := N(K)\cap G(K)^+$. Then $N(K)^+$ maps onto the Weyl group $W :=N(K)/Z(K)$ of $G$ (this can be seen using, for example, [CGP, Prop.\,C.2.24(i)]).   
\vskip1mm

 Let $\cA$ be the apartment of $\sB$ corresponding to $S$. As in 2.14, let $\varphi: N(K)\rightarrow {\rm{Aff}}(\cA)$ be the action map, then the affine Weyl group $W_{\rm{aff}}$ of $G/K$ is by definition the subgroup $\varphi(N(K)^+)$ of ${\rm{Aff}}(\cA)$. 
 

\vskip3mm

\ni{\bf 3. Bruhat-Tits theory for $G$ over $K$}  
\vskip3mm

\ni {\bf 3.1.}  As before, we will continue to assume in this section that $H$ is semi-simple and the central torus of $G$ is $K$-anisotropic (see 2.13). For a $\Theta$-stable nonempty bounded subset $\Omega$ of an apartment of $\cB(H/K)$, let $\sG^{\circ}_{\Omega}$ be the Bruhat-Tits group scheme as defined in 2.5(i).  We will denote the special fiber of $\sG^{\circ}_{\Omega}$ by $\osG^{\circ}_{\Omega}$ in the sequel. 
\vskip1mm

Given a point  $x\in \sB$, for simplicity we will denote $\sG^{\circ}_{\{x\}}$, $\sH_{\{x\}}$, $\sH^{\circ}_{\{x\}}$ and $\sH^{\Theta}_{\{x\}}$ by $\sG^{\circ}_x$, $\sH_x$, $\sH^{\circ}_x$ and $\sH^{\Theta}_x$ respectively, and the special fibers of these group schemes will be denoted by ${\overline\sG}^{\circ}_x$, ${\overline\sH}_x$, ${\overline\sH}^{\circ}_x$ and ${\overline\sH}^{\Theta}_x$ respectively. The subgroup of $H(K)$ (resp.\,$G(K)$) consisting of elements that fix $x$ will be denoted by $H(K)^x$ (resp.\,$G(K)^x$). The subgroup $\sG^{\circ}_x(\cO)\,(\subset G(K)^x)$ is of finite index in $G(K)^x$.  
\vskip2mm

\ni{\bf 3.2.} Let $\Omega'\prec \Omega$ be nonempty $\Theta$-stable bounded subsets of an apartment of $\cB(H/K)$. The $\cO$-group scheme homomorphism $ \sH_{\Omega}\rightarrow \sH_{\Omega'}$ of [P2,\,1.10]  restricts to a homomorphism $\rho_{\Omega', \Omega}: \sH^{\circ}_{\Omega}\rightarrow \sH^{\circ}_{\Omega'}$, and by [E, Prop.\,3.5], or [CGP, Prop.\,A.8.10(2)], it induces a $\cO$-group scheme homomorphism $\sH_{\Omega}^{\Theta}\rightarrow \sH_{\Omega'}^{\Theta}$. The last homomorphism gives a $\cO$-group scheme homomorphism $\rho^G_{\Omega',\Omega}: (\sH^{\Theta}_{\Omega})^{\circ}=\sG^{\circ}_{\Omega}\rightarrow  
\sG^{\circ}_{\Omega'}=(\sH_{\Omega'}^{\Theta})^{\circ}$ that is the identity homomorphism on the generic fiber $G$.  
\vskip2mm

\ni{\bf 3.3.}  Let $\cA$ be the apartment of $\sB$ corresponding to a maximal $K$-split torus $S$ of $G$ and $\Omega$ be a nonempty bounded subset of  $\cA$. The apartment $\cA$ is contained in an apartment $A$ of $\cB(H/K)$ that corresponds to a maximal $K$-split torus $T$ of $H$ containing $S$ and $\cA = A\cap \sB =A^{\Theta}$\,(2.16).  So $\Omega$ is a bounded subset of $A$.   The group scheme $\sH_{\Omega}$ contains a closed split $\cO$-torus $\sT$ with generic fiber $T$, see [P2, 1.9].  Let $\sS$ be the closed $\cO$-subtorus of $\sT$ whose generic fiber is $S$ ($\sS$ is the schematic closure of $S$ in $\sT$).   The automorphism group $\Theta$ of $\sH_{\Omega}$ acts trivially on the  $\cO$-torus $\sS$ (since $S\subset G\subset H^{\Theta}$) and hence $\sS$ is contained in $\sG^{\circ}_{\Omega}$. The special fiber $\osS$ of $\sS$ is a maximal torus of $\osG^{\circ}_{\Omega}$ since $S$ is a maximal $K$-split torus of $G$.  
\vskip1mm

For $x\in \sB$, let $P_x := \sG^{\circ}_x(\cO)$ be the parahoric subgroup of $G(K)$ associated with the point $x$. Let $S$ be a maximal $K$-split torus of $G$ such that $x$ lies in the apartment $\cA$ of $\sB$ corresponding to $S$. Then the group scheme $\sG^{\circ}_x$ contains a closed split $\cO$-torus $\sS$ whose generic fiber is $S$. 
  
  \vskip2mm

\ni{\bf Proposition 3.4.} {\em Let $\cA$ and $\cA'$ be apartments of $\sB$ and $\Omega$ a nonempty bounded subset of $\cA\cap \cA'$. Then there exists an element  $g\in \sG^{\circ}_{\Omega}(\cO)$ that maps $\cA$ onto $\cA'$. Any such element fixes $\Omega$ pointwise.}
\vskip2mm

\ni{\it Proof.} We will use Proposition 2.1(ii)  of [P2], with $\cO$ in place of $\fo$,  and denote 
$\sG^{\circ}_{\Omega}$ by $\sG$, and its special fiber by $\osG$, in this proof. Let $S$ and $S'$ be the maximal $K$-split tori of $G$ corresponding to the apartments $\cA$ and $\cA'$ respectively and $\sS$ and $\sS'$ be the closed $\cO$-tori of $\sG$ with generic fibers $S$ and $S'$  respectively.  The special fibers $\osS$ and $\osS'$ of $\sS$ and $\sS'$ are maximal split tori of $\osG$, and hence  according to a result of Borel and Tits there is an element $\overline{g}$ of $\osG(\kappa)$ which conjugates $\osS$ onto $\osS'$ [CGP,\:Thm.\,C.2.3].  Now [P2,\:Prop.\,2.1(ii)] implies that there exists a $g\in\sG(\cO)$ lying over $\overline{g}$ that conjugates $\sS$ onto $\sS'$. This element fixes $\Omega$ pointwise and conjugates $S$ onto $S'$ and hence maps  $\cA$ onto $\cA'$.  \hfill$\Box$
\vskip1.5mm

\ni{\bf 3.5. Polysimplicial structure on $\sB$.} Let $P$ be a parahoric subgroup of $G(K)$ and $\sG^{\circ}$ be the Bruhat-Tits parahoric $\cO$-group scheme associated with $P$\,(2.5(ii)). Let $\cB(H/K)^P$ denote the set of points of $\cB(H/K)$ fixed by $P$. According to Corollary 2.3, $\cB(H/K)^P$ is the union of facets pointwise fixed by  $P$. Let $\overline{\cF}_P:= \cB(H/K)^P\cap \sB$. This closed convex subset  is by definition the {\it closed facet} of $\sB$ associated with the parahoric subgroup $P$.  The  $\cO$-group scheme $\sG^{\circ}$ contains a closed split $\cO$-torus $\sS$ whose generic fiber $S$ is a maximal $K$-split torus of $G$\,(3.3). The subgroup  $\sS(\cO)$ (of $S(K)$) is the maximal bounded subgroup of $S(K)$ [PY2, 3.6] and it is contained in $P\,(=\sG^{\circ}(\cO))$,  so,  according to Corollary 2.3,  $\cB(H/K)^P$ is contained in the enlarged building  $\cB(Z_H(S)/K)$ of $Z_H(S)(K)$. This implies that  the closed facet $\overline{\cF}_P$ is contained in the apartment $\cA:=\cB(Z_H(S)/K)^{\Theta}\,(=\cB(Z_H(S)/K)\cap \sB)$  of $\sB$ corresponding to the maximal $K$-split torus $S$ of $G$.   
\vskip1mm

 Given another parahoric subgroup  subgroup $Q$ of $G(K)$, if ${\overline\cF}_Q = {\overline\cF}_P$,  then $Q = P$. (To see this, we choose points $x,\,y\in \sB$ such that $\sG^{\circ}_x(\cO) = P$ and $\sG^{\circ}_y (\cO)= Q$. Then $y\in \overline{\cF}_Q =\overline{\cF}_P$. So $P$ fixes $y$. Now using 2.5(iv) we see that $P\subset Q$. We  similarly see  that $Q\subset P$.) Hence if $Q\supsetneq P$, then  ${\overline\cF}_Q$ is properly contained in ${\overline\cF}_P$. Let $\cF_P$ be the subset of points of $\overline{\cF}_P$ that are not fixed by any parahoric subgroup of $G(K)$ larger than $P$. Then $\cF_P  = \overline{\cF}_P - \bigcup_{Q\supsetneq P}\overline{\cF}_Q$.   By definition, $\cF_P$ is the {\it facet} of $\sB$ associated with the parahoric subgroup $P$ of $G(K)$, and as $P$ varies over the set of parahoric subgroups of $G(K)$, these are are all the facets of $\sB$. We will show below  (Propositions 3.8 and 3.10) that $\cF_P$ is  convex and bounded. 

\vskip2mm

In the following two lemmas (3.6 and 3.7),  $\kappa$ is any field of characteristic $p\geqslant 0$. We  will use the notation introduced in [CGP, \S2.1]. 

\vskip2mm

\ni{\bf Lemma 3.6.} {\em Let $\cH$ be a smooth connected affine algebraic $\kappa$-group and $\cQ$ be a pseudo-parabolic $\kappa$-subgroup of $\cH$. Let $\cS$ be a $\kappa$-torus of $\cQ$ whose image in the maximal pseudo-reductive quotient $\cM:=\cQ/\sR_{u,\kappa}(\cQ)$ of $\cQ$ contains the maximal central torus of $\cM$. Then any $1$-parameter subgroup $\lambda: {\rm{GL}}_1\rightarrow \cH$ such that $\cQ = P_{\cH}(\lambda)\sR_{u,\kappa}(\cH)$ has a conjugate under $\sR_{u,\kappa}(\cQ)(\kappa)$ with image in $\cS$.}
\vskip2mm

\ni{\it Proof.}  Let $\lambda: {\rm{GL}}_1\rightarrow \cH$ be a $1$-parameter subgroup such that $\cQ = P_{\cH}(\lambda)\sR_{u,\kappa}(\cH)$. The image $\cT$ of $\lambda$ is contained in $\cQ$ and it maps into the central torus of  $\cM$. Therefore, $\cT$ is contained in the solvable subgroup $\cS\sR_{u, \kappa}(\cQ)$ of $\cQ$. Note that as $\cS$ is commutative, the derived subgroup of $\cS\sR_{u, \kappa}(\cQ)$ is contained in $\sR_{u, \kappa}(\cQ)$, so  the maximal $\kappa$-tori of $\cS\sR_{u, \kappa}(\cQ)$ are conjugate to each other under $\sR_{u, \kappa}(\cQ)(\kappa)$\,[Bo,\:Thm.\,19.2].  Hence,  there is a $u\in\sR_{u,\kappa}(\cQ)(\kappa)$ such that $u\cT u^{-1}\subset \cS$. Then the image of the $1$-parameter subgroup $\mu: {\rm{GL}}_1\rightarrow \cS$, defined as $\mu(t) =u\lambda(t)u^{-1}$, is contained in $\cS$. \hfill$\Box$
 \vskip2mm

\ni{\bf Lemma 3.7.} {\em Let $\cH$ be a smooth connected affine algebraic $\kappa$-group given with an action by a finite group $\Theta$. We assume that $p$ does not divide the order of $\Theta$. Let $\cG= {(\cH^{\Theta})}^{\circ}$. Then 
\vskip1mm

{\rm{(i)}} $\sR_{u, \kappa}(\cG)= (\cG\cap \sR_{u, \kappa}(\cH))^{\circ}= {(\sR_{u, \kappa}(\cH)^{\Theta})}^{\circ}$; moreover,  $\cG/(\cG\cap \sR_{u,\kappa}(\cH))$ is pseudo-reductive, and if  $\kappa$ is perfect then $\cG\cap \sR_{u, \kappa}(\cH) =\sR_{u,\kappa}(\cG)$. 

\vskip.7mm 
{\rm{(ii)}} Given a $\Theta$-stable  pseudo-parabolic $\kappa$-subgroup $\cQ$ of $\cH$, $\cP:=\cG\cap \cQ$ is a pseudo-parabolic $\kappa$-subgroup of \,$\cG$, so $\cP$ is connected and it equals ${(\cQ^{\Theta})}^{\circ}$. 

\vskip.7mm
{\rm{(iii)}} Conversely, given a pseudo-parabolic $\kappa$-subgroup $\cP$ of $\cG$, and a maximal $\kappa$-torus $\cS\subset \cP$, there is a $\Theta$-stable pseudo-parabolic $\kappa$-subgroup $\cQ$ of $\cH$, $\cQ$ containing the centralizer $Z_{\cH}(\cS)$ of $\cS$ in $\cH$, such that $\cP = \cG\cap \cQ = {(\cQ^{\Theta})}^{\circ}$.}
\vskip2mm

\ni{\it Proof.}  The first assertion of (i) immediately follows from [CGP, Prop.\,A.8.14(2)]. Now we observe that as $\sR_{u,\kappa}(\cG) = (\cG\cap \sR_{u,\kappa}(\cH))^{\circ}$, the quotient $(\cG\cap \sR_{u,\kappa}(\cH))/\sR_{u,\kappa}(\cG)$ is a finite \'etale (unipotent) normal subgroup of 
$\cG/\sR_{u,\kappa}(\cG)$ so it is central. Thus the kernel of  the quotient map $\pi: \cG/\sR_{u,\kappa}(\cG)\rightarrow \cG/(\cG\cap \sR_{u,\kappa}(\cH))$ is an \'etale unipotent central subgroup. Hence, $ \cG/(\cG\cap \sR_{u,\kappa}(\cH))$ is pseudo-reductive as $\cG/\sR_{u,\kappa}(\cG)$ is. Moreover, if $\kappa$ is perfect then every pseudo-reductive $\kappa$-group is reductive and such a group does not contain a nontrivial  \'etale unipotent normal subgroup. This implies that if $\kappa$ is perfect, then $\sR_{u, \kappa}(\cG)= \cG\cap \sR_{u, \kappa}(\cH)$. 
\vskip1mm

Since $\sR_{u,\kappa}(\cG)\subset \cG\cap \sR_{u,\kappa}(\cH)\subset \cG\cap \cQ$, to prove (ii), we can replace $\cH$ by its  pseudo-reductive quotient $\cH/\sR_{u, \kappa}(\cH)$ and assume that $\cH$ is pseudo-reductive. Then $\cG$ is also pseudo-reductive by (i). The $\kappa$-unipotent radical  $\sR_{u, \kappa}(\cQ)$  of $\cQ$ is $\Theta$-stable. Let $\cS$ be a $\Theta$-stable $\kappa$-torus in $\cQ$ that maps isomorphically onto the maximal central torus of the pseudo-reductive quotient $\overline{\cQ}:= \cQ/\sR_{u, \kappa}(\cQ)$ (Lemma 2.7). By Lemma 3.6, there exists a 1-parameter subgroup $\lambda: {\rm{GL}}_1\rightarrow \cS$ such that $\cQ = P_{\cH}(\lambda)$. Let $\mu =\sum_{\theta\in \Theta} \theta\cdot \lambda$. Then $\mu$ is invariant under $\Theta$ and so it is a 1-parameter subgroup of $\cG$. We will now show that $\cQ = P_{\cH}(\mu)$. Let $\Phi$ (resp.\,$\Psi$) be the set of weights in the Lie algebra of $\cQ$ (resp.\,$P_{\cH}(\mu)$) with respect to the adjoint action of $\cS$. Then since $\cQ$, $P_{\cH}(\mu)$  and $\cS$ are $\Theta$-stable, the subsets $\Phi$ and $\Psi$ (of ${\rm{X}}(\cS)$) are stable under the action of $\Theta$ on ${\rm{X}}(\cS)$. Hence, for all $a \in \Phi$, as $\langle a,\lambda\rangle \geqslant 0$, we conclude that  $\langle a,\mu\rangle \geqslant 0$. Therefore, $\Phi\subset \Psi$.  On the other hand, for $b\in \Psi$, $\langle b,\mu\rangle \geqslant 0$. If $b\,(\in \Psi)$ does not belong to $\Phi$, then for $\theta\in \Theta$, $\theta\cdot b\notin \Phi$, so for all $\theta \in \Theta$, $\langle \theta\cdot b,\lambda\rangle <0$, which implies that $\langle b,\mu\rangle < 0$. This is a contradiction. Therefore, $\Phi = \Psi$ and so $\cQ = P_{\cH}(\mu)$.  Now observe that ${(\cQ^{\Theta})}^{\circ}\subset \cG\cap\cQ\subset \cQ^{\Theta}$. As $\cQ^{\Theta}$ is a smooth subgroup ([E, Prop.\,3.4] or [CGP, Prop.\,A.8.10(2)]), $\cG\cap \cQ$ is a smooth $\kappa$-subgroup, and since it contains the pseudo-parabolic $\kappa$-subgroup $P_{\cG}(\mu)$ of $\cG$, it is a pseudo-parabolic $\kappa$-subgroup of $\cG$ [CGP, Prop.\,3.5.8], hence in particular it is connected.     Therefore, $\cG\cap \cQ = {(\cQ^{\Theta})}^{\circ}$.   
\vskip1mm

Now we will prove (iii). Let $\lambda: {\rm{GL}}_1\rightarrow \cS$ be a $1$-parameter subgroup such that $\cP = P_{\cG}(\lambda)\sR_{u, \kappa}(\cG)$. Then $\cQ : = P_{\cH}(\lambda)\sR_{u, \kappa}(\cH)$  is a pseudo-parabolic $\kappa$-subgroup of $\cH$ that is $\Theta$-stable (since $\lambda$ is $\Theta$-invariant) and it contains $\cP$ as well as $Z_{\cH}(\cS)$. According to (ii), $\cG\cap \cQ = {(\cQ^{\Theta})}^{\circ}$ is a pseudo-parabolic $\kappa$-subgroup of $\cG$ containing $\cP$. The Lie algebras of $\cP$ and ${(\cQ^{\Theta})}^{\circ}$ are clearly equal. This implies that $\cP = \cG\cap \cQ={(\cQ^{\Theta})}^{\circ}$ and we have proved (iii). \hfill$\Box$

\vskip2mm
\ni{\bf Proposition 3.8.} {\em Let $P$ be a parahoric subgroup of $G(K)$ and $\cF_P$ and $\overline{\cF}_P$  be as in $3.5$.}
\vskip.5mm

 (i) {\em Given $x\in \cF_P$ and $y\in \overline{\cF}_P$, for every point $z$ of the geodesic $[xy]$, except possibly for $z=y$, $\sG^{\circ}_z(\cO) =P$.}
\vskip.5mm

(ii) {\em Let $F$ be a facet of $\cB(H/K)$ that meets $\overline{\cF}_P$ and is maximal among such facets. Then $\sG^{\circ}_F(\cO) =P$. Thus $F\cap \sB\subset \cF_P$.}
\vskip2mm

The first assertion of this proposition implies that $\cF_P$ is convex.  The second assertion implies that $\cF_P$ is an open-dense subset of $\overline{\cF}_P$, hence the closure of $\cF_P$ is $\overline{\cF}_P$.
\vskip2mm

\ni{\it Proof.} To prove the first assertion, let $[xy]$ be the geodesic joining $x$ and $y$. Let $F_0,\, F_1, \ldots,\, F_n$ be the facets of $\cB(H/K)$ containing a segment of positive length of the geodesic $[xy]$  (so each $F_i$ is $\Theta$-stable and  is fixed pointwise by $P$, hence $P\subset \sG^{\circ}_{F_i}(\cO)$, cf.\,2.5(iv)). Then $[xy]\subset \bigcup_i \overline{F}_i$. We assume  the facets $\{F_i\}$ indexed so that $x$ lies in $\overline{F}_0$, $y$ lies in $\overline{F}_n$,  and  for  each $i<n$, $\overline{F}_i \cap \overline{F}_{i+1}$ is nonempty.   Let $z_0 = x$. For every positive integer $i\,(\leqslant n)$, $\overline{F}_{i-1} \cap \overline{F}_{i}$  contains a unique point of $[xy]$; we will denote this point by $z_{i}$. 
\vskip1mm

To prove the second assertion of the proposition along with the first, we take $x$ to be a point of $\sB$ such that $\sG^{\circ}_x (\cO) = P$ (so $x\in \cF_P$) and take $y$ to be any point of $F\cap \sB$. Let $[xy]$, and for $i\leqslant n$, $F_i$ and $z_i$  be as in the preceding paragraph. Then $F_n = F$. 
\vskip1mm

Since $x\in {\overline{F}}_0$, there is a $\cO$-group scheme homomorphism $\sG^{\circ}_{F_0}\rightarrow\sG^{\circ}_x$ that is the identity on the generic fiber $G$. 
Thus, $\sG^{\circ}_{F_0}(\cO)\subset P$.  But $P\subset \sG^{\circ}_{F_0}(\cO)$, so $\sG^{\circ}_{z_0}(\cO)=\sG^{\circ}_{F_0}(\cO) = P$. Let  $j\,(\leqslant n)$ be a positive  integer such that for all $i< j$, $\sG^{\circ}_{z_i}(\cO)=\sG^{\circ}_{F_i}(\cO) = P$. The inclusion of $\{z_{j}\}$ in $\overline{F}_{j-1} \cap \overline{F}_{j}$ gives rise to $\cO$-group scheme homomorphisms $\sH_{F_{j-1}}\stackrel{\sigma_{j}}\longrightarrow \sH_{z_{j}}\stackrel{\rho_{j}}\longleftarrow \sH_{F_{j}}$ that are the identity on the generic fiber $H$. The images of the induced homomorphisms $\osH^{\circ}_{F_{j-1}}\stackrel{{\overline\sigma}_{j}}\longrightarrow \osH^{\circ}_{z_{j}}\stackrel{{\overline\rho}_{j}}\longleftarrow \osH^{\circ}_{F_{j}}$ are pseudo-parabolic $\kappa_s$-subgroups of $\osH^{\circ}_{z_{j}}$\,([P2,\,1.10(2)]). We conclude by Lie algebra consideration that ${\overline{\sigma}}_{j}(\osG^{\circ}_{F_{j-1}})={({\overline{\sigma}}_{j}(\osH^{\circ}_{F_{j-1}})^{\Theta})}^{\circ}$ and ${\overline{\rho}}_{j}(\osG^{\circ}_{F_{j}})={({\overline{\rho}}_{j}(\osH^{\circ}_{F_{j}})^{\Theta})}^{\circ}$, and Lemma 3.7(ii) implies that both of these subgroups are pseudo-parabolic subgroups of $\osG^{\circ}_{z_{j}}$. As $\sG^{\circ}_{F_{j-1}}(\cO) = P$, whereas,  $P\subset \sG^{\circ}_{F_{j}}(\cO)\, (\subset  \sG^{\circ}_{z_{j}}(\cO))$, we see that ${\overline{\sigma}}_{j}(\osG^{\circ}_{F_{j-1}})$ is contained in ${\overline{\rho}}_{j}(\osG^{\circ}_{F_{j}})$.  Let $\overline Q$ and $\overline{Q}'$ respectively be the images of \,${\overline{\sigma}}_{j}(\osG^{\circ}_{F_{j-1}})$ and ${\overline{\rho}}_{j}(\osG^{\circ}_{F_{j}})$  in  the maximal pseudo-reductive quotient  $\orG_{z_{j}} := \osG^{\circ}_{z_{j}}/\sR_{u,\kappa_s}(\osG^{\circ}_{z_{j}})$ of $\osG^{\circ}_{z_{j}}$. Then $\overline{Q}\subset \overline{Q}'$, and both of them are pseudo-parabolic subgroups of $\orG_{z_{j}}$.  
\vskip1mm

Now let $S$ be a maximal $K$-split torus of $G$ such that the apartment of $\sB$ corresponding to $S$ contains the geodesic $[xy]$ and let  $v\in V(S)$ so that $v+x = y$. Then for all sufficiently small positive real number $\epsilon$, $-\epsilon v+z_{j}\in F_{j-1}$ and $\epsilon v+z_{j}\in F_{j}$. Using  [P2,\,1.10(3)] we infer that the images of the pseudo-parabolic subgroups ${\overline{\sigma}}_j(\osH^{\circ}_{F_{j-1}})$ and ${\overline{\rho}}_{j}(\osH^{\circ}_{F_{j}})$ (of $\osH^{\circ}_{z_j}$) in the maximal pseudo-reductive quotient $\opredH_{z_{j}}  :=  \osH^{\circ}_{z_{j}}/\sR_{u,\kappa_s}(\osH^{\circ}_{z_{j}})$ of $\osH^{\circ}_{z_{j}}$ are opposite pseudo-parabolic subgroups. Therefore, the image $\cH$ of \,${\overline{\sigma}}_j(\osH^{\circ}_{F_{j-1}})\cap {\overline{\rho}}_{j}(\osH^{\circ}_{F_{j}})$ in $\opredH_{z_{j}}$ is pseudo-reductive. Proposition A.8.14\,(2) of [CGP] implies then that ${(\cH^{\Theta})}^{\circ}$ is pseudo-reductive.  It is obvious that under the natural homomorphism 
$\pi: \orG_{z_j}\rightarrow  \opredH_{z_j}$, the image of $\overline{Q}= \overline{Q}\cap \overline{Q}'$ is ${(\cH^{\Theta})}^{\circ}$.  As the kernel of the homomorphism $\pi$  is a finite  (\'etale unipotent) subgroup (Lemma 3.7(i)), and ${(\cH^{\Theta})}^{\circ}$ is pseudo-reductive, we see that $\overline{Q}$ is a pseudo-reductive subgroup of $\orG_{z_{j}}$. But since $\overline{Q}$ is a pseudo-parabolic subgroup of the latter, we must have $\overline{Q} = \orG_{z_{j}}$, and hence, $\overline Q' = \orG_{z_{j}}$.  So, 
${{\overline\sigma}}_j(\osG^{\circ}_{F_{j-1}}) = \osG^{\circ}_{z_{j}}= {{\overline\rho}}_{j}(\osG^{\circ}_{F_{j}})$. 
\vskip1mm

Since the natural homomorphism $\sG^{\circ}_{F_{j-1}}(\cO)\rightarrow \osG^{\circ}_{F_{j-1}}(\kappa)$ is surjective (as $\cO$ is henselian and $\sG^{\circ}_{F_{j-1}}$ is smooth,\,[EGA\,IV$_4$, 18.5.17]),  and ${{\overline\sigma}}_j(\osG^{\circ}_{F_{j-1}}) = \osG^{\circ}_{z_{j}}$, the image of $\sG^{\circ}_{F_{j-1}}(\cO)\,(\subset \sG^{\circ}_{z_j}(\cO))$ in $\osG^{\circ}_{z_j}(\kappa)$ is Zariski-dense in $\osG^{\circ}_{z_{j}}$. From this we see that $$\cO[\sG^{\circ}_{z_j}] =\{f\in K[G]\,\vert\, f(\sG^{\circ}_{F_{j-1}}(\cO))\subset\cO\} = \cO[\sG^{\circ}_{F_{j-1}}],$$ cf.\,[BrT2,\,1.7.2] and 2.1. Therefore,  $\sigma_j\vert_{\sG^{\circ}_{F_{j-1}}}:\sG^{\circ}_{F_{j-1}}\rightarrow \sG^{\circ}_{z_j}$ is a $\cO$-group scheme isomorphism. We similarly see that  $\rho_j\vert_{\sG^{\circ}_{F_{j}}}:\sG^{\circ}_{F_{j}}\rightarrow \sG^{\circ}_{z_j}$ is a $\cO$-group scheme isomorphism.  Now since $\sG^{\circ}_{F_{j-1}}(\cO) = P$, we conclude that $P=\sG^{\circ}_{z_{j}}(\cO) = \sG^{\circ}_{F_{j}}(\cO)$. By  induction it follows that $P=\sG^{\circ}_{z_i}(\cO)=\sG^{\circ}_{F_{i}}(\cO)$ for all $i\leqslant n$. In particular, for all $z\in [xy]$, except possibly for $z=y$, $\sG^{\circ}_z(\cO) = P$, and $\sG^{\circ}_{F_n}(\cO) = P$. \hfill$\Box$

\vskip2mm

For parahoric subgroups $P$ and $Q$ of $G(K)$, if $\cF_P\cap\cF_Q$ is nonempty, then for any $z$ in this intersection, $P = \sG^{\circ}_z(\cO)= Q$ (Proposition 3.8(i)).  Thus every point of $\sB$ is contained in a unique facet.
\vskip1mm

For a parahoric subgroup $Q$ of $G(K)$ containing $P$, obviously, $\cF_Q\subset \overline{\cF}_Q\subset \overline{\cF}_P$, thus $\cF_Q\prec \cF_P$ and hence $\cF_P$ is a maximal facet if and only if $P$ is a minimal parahoric subgroup of $G(K)$. The maximal facets of $\sB$ are called the {\it chambers} of $\sB$. It is easily seen using the observations contained in 2.9 that all the chambers are of equal dimension. 
We say that a facet $\cF'$ of $\sB$ is a {\it face} of a facet $\cF$ if \,$\cF'\prec \cF$, i.e., if \,$\cF'$ is contained in the closure of \,$\cF$.

\vskip2mm

We will use the following simple lemma in the proof of the next proposition.

\vskip2mm

\ni{\bf Lemma 3.9.} {\em Let $S$ be a maximal $K$-split torus of $G$, $\cA$ the corresponding apartment of $\sB$, and $\cC$ be a noncompact closed convex subset of $\cA$. Then for any point $x\in \cC$, there is an infinite ray originating at $x$ and contained in $\cC$.} 
\vskip2mm

\ni{\it Proof.} Recall that $\cA$ is an affine space under the vector space $V(S) = \bR\otimes_{\bZ}{\rm{X}}_{\ast}(S)$. We identify $\cA$ with $V(S)$ using translations by elements in the latter, with $x$ identified with the origin $0$, and use a positive definite inner product on $V(S)$ to get a norm on $\cA$. With this identification, $\cC$ is a closed convex subset of $V(S)$ containing $0$. Since $\cC$ is noncompact, there exist unit vectors $v_i\in V(S)$, $i\geqslant 1$, and positive real numbers $s_i\rightarrow \infty$ such that $s_iv_i$ lies in $\cC$. After replacing $\{v_i\}$ by a subsequence,  we may (and do) assume that the sequence  $\{v_i\}$  converges to a unit vector $v$. We will now show that for every nonnegative real number $t$, $tv$ lies in $\cC$, this will prove the lemma. To see that $tv$ lies in $\cC$, it suffices to observe that for a given $t$, the sequence $\{tv_i\}$  converges to $tv$, and for all sufficiently large $i$ (so that $s_i\geqslant t$), $tv_i$ lies in $\cC$. \hfill$\Box$
\vskip2mm

 \ni{\bf Proposition  3.10.} {\em For any parahoric subgroup $P$ of $G(K)$, the associated closed facet ${\overline\cF}_P$ of $\sB$, and so also the associated facet $\cF_P\,(\subset {\overline\cF}_P)$, is bounded.}
 \vskip2mm
 
 \ni{\it Proof.}  Let $S$ be a maximal $K$-split torus of $G$ such that the corresponding apartment of $\sB$ contains $\overline{\cF}_P$\,(3.5). Assume, if possible, that ${\overline\cF}_P$ is noncompact and fix a point $x$ of $\cF_P$. Then, according to the preceding lemma, there is an infinite ray $\cR :=\{tv+x\, |\, t\in \bR_{\geqslant 0}\}$, for some $v\in V(S)$, originating at $x$ and contained in $\overline{\cF}_P$. It is obvious from Proposition 3.8(i) that this ray is actually contained in $\cF_P$. Hence, for every point $z\in \cR$, $\sG^{\circ}_z(\cO) = P$.  
 \vskip.5mm
 
 As the central torus of $G$ has been assumed to be $K$-anisotropic, there is a nondivisible root $a$ of $G$, with respect to $S$, such that $\langle a,\, v\rangle>0$.  Let $S_a$ be the identity component of the kernel of $a$ and $G_a$ (resp.\,$H_a$) be the derived subgroup of the centralizer of $S_a$ in $G$ (resp.\,$H$).   Fix $t\in \bR_{\geqslant 0}$, and let  $y= tv+x\in \cR$.  Let $\sS$ be the closed $1$-dimensional $\cO$-split torus of $\sG^{\circ}_y$ whose generic fiber is the maximal $K$-split torus of $G_a$ contained in $S$ and let  $\lambda: {\rm{GL}}_1\rightarrow \sS\,(\hookrightarrow \sG^{\circ}_y\hookrightarrow \sH_y)$ be the $\cO$-isomorphism such that $\langle a,\lambda\rangle>0$. Let $c= \langle a, v\rangle/\langle a, \lambda\rangle$. Then $\langle a, v-c\lambda\rangle = 0$.  
 \vskip.5mm
 
Let $\sU_y$ be the $\cO$-subgroup scheme of $\sH_y$ representing the functor  $$R \leadsto \{ h\in \sH_y(R)\,|\, \lim_{t \to 0} \lambda(t)h\lambda(t)^{-1} = 1\},$$ cf.\,[CGP, Lemma 2.1.5].  Using the last assertion of 2.1.8(3), and the first assertion of 2.1.8(4), of [CGP] (with $k$, which is an an arbitrary commutative ring in these assertions, replaced by $\cO$, and $G$ replaced by $\sH_y$), we see that  $\sU_y$ is a closed smooth unipotent $\cO$-subgroup scheme of $\sH_y$ with connected fibers; the generic fiber of $\sU_y$  is $U_H(\lambda)$, where $U_H(\lambda)$ is as in [CGP, Lemma 2.1.5] with $G$ replaced by $H$. We consider the smooth closed $\cO$-subgroup scheme $\sU_y^{\Theta}$ of $\sU_y$. As $\sU_y^{\Theta}$ is clearly normalized by $\sS$, it has connected fibers, and hence it is contained in $(\sH_y^{\Theta})^{\circ}= \sG_y^{\circ}$.  The generic fiber of $\sU_y^{\Theta}$ is $U_H(\lambda)^{\Theta}$ that contains the root group $U_a\,(=U_{G_a}(\lambda))$ of $G$ corresponding to the root $a$.  

As $\bigcup_{z\in \cR} \sU_z(\cO) \supset U_{H_a}(\lambda)(K)\supset U_a(K)$, we see that $\bigcup_{z\in \cR}\sU^{\Theta}_z(\cO)\supset U_a(K)$. Now since  $\sG^{\circ}_z\supset \sU^{\Theta}_z$, we conclude that $\bigcup_{z\in \cR}\sG^{\circ}_z(\cO)\supset U_a(K)$. But for all $z\in \cR$, $\sG^{\circ}_z(\cO) = P$, so the parahoric subgroup $P$ contains the unbounded subgroup $U_a(K)$. This is a contradiction.\hfill$\Box$ 
\vskip2mm

Proposition 3.10 implies that each closed facet of $\sB$ is a compact polyhedron. Considering the facets lying on the boundary of a maximal closed facet of $\sB$, we see that $\sB$ contains facets of every dimension $\leqslant K$-rank\,$G$. 
\vskip2mm

\ni{\bf 3.11.} Let $P$ be a parahoric subgroup of $G(K)$ and $\cF:=\cF_P$ be the facet of $\sB$ associated to $P$ in 3.5. Then for any $x\in \cF$, since $P\subset \sG^{\circ}_{\cF} (\cO)\subset \sG^{\circ}_x(\cO)= P$\, (2.5(iv)), $\sG^{\circ}_{\cF}(\cO) = P$ and hence the natural  $\cO$-group scheme homomorphism $\sG^{\circ}_{\cF}\rightarrow \sG^{\circ}_x$ is an isomorphism. In particular, for any facet $F$ of $\cB(H/K)$ that meets $\cF$, $\sG^{\circ}_{\cF} = \sG^{\circ}_F$. 

\vskip2mm
\ni{\bf Proposition 3.12.} {\em Let $\cF$ be a facet of $\sB$. Then the $\kappa$-unipotent radical  $\sR_{u, \kappa}(\osG^{\circ}_{\cF})$ of $\osG^{\circ}_{\cF}$ equals ${(\osG^{\circ}_{\cF}\cap\sR_{u,\kappa}(\osH^{\circ}_{\cF}))}^{\circ}$.}
\vskip1.5mm

{\em Let $\cF$ and $\cF'$ be two  facets of $\sB$, with $\cF'\prec \cF$.  Then$:$ 
\vskip1mm

{\rm{(i)}} The kernel of the induced  homomorphism  ${\overline{\rho}}^{G}_{\cF',\cF}: \osG^{\circ}_{\cF}\rightarrow \osG^{\circ}_{\cF'}$ between the special fibers  is a smooth unipotent $\kappa$-subgroup of \,$\osG^{\circ}_{\cF}$ and the image $\fp(\cF'/\cF)$ is a pseudo-parabolic $\kappa$-subgroup of  $\osG^{\circ}_{\cF'}$. 
\vskip1mm

{\rm{(ii)}} If $F$ and $F'$ are facets of $\cB(H/K)$, $F'\prec F$, that meet $\cF$ and $\cF'$ respectively, then $\fp(\cF'/\cF) = (\osQ^{\Theta})^{\circ}$, where $\osQ$ is the image of ${\overline{\rho}}_{F',F}:\osH^{\circ}_{F}\rightarrow \osH^{\circ}_{F'}$. 
\vskip1mm

{\rm{(iii)}} The inverse image of the subgroup $\fp(\cF'/\cF)(\kappa)$ of \,$\osG^{\circ}_{\cF'}(\kappa)$, under the natural surjective homomorphism $\sG^{\circ}_{\cF'} (\cO)\rightarrow \osG^{\circ}_{\cF'}(\kappa)$, is  $\rho^G_{\cF',\cF}(\sG^{\circ}_{\cF}(\cO))\,(\subset \sG^{\circ}_{\cF'}(\cO))$.
\vskip1mm

Given a pseudo-parabolic $\kappa$-subgroup $\osP$ of $\osG^{\circ}_{\cF'}$, there is a facet $\cF$ of $\sB$ with $\cF'\prec \cF$ such that the image of the  homomorphism  $\overline{\rho}^G_{\cF',\cF}: \osG^{\circ}_{\cF}\rightarrow   \osG^{\circ}_{\cF'}$ equals $\osP$.}
\vskip2mm

 \ni{\it Proof.} The first assertion of the proposition follows immediately from Lemma 3.7(i). 
 \vskip1mm
 
 To prove  (i), we fix $x\in \cF'$ and let $F'$ be the facet of $\cB(H/K)$ containing $x$. As the closure of $\cF$ contains $x$, there is a facet $F$ of $\cB(H/K)$ that meets $\cF$ and contains $x$ in its closure. Then  $F' \subset \overline{F}$, i.e., $F'\prec F$, and $F$ and $F'$ meet $\cF$ and $\cF'$ respectively. Hence, $\sG^{\circ}_{\cF} = \sG^{\circ}_F= {(\sH_F^{\Theta})}^{\circ}$ and $\sG^{\circ}_{\cF'} = \sG^{\circ}_{F'}= {(\sH_{F'}^{\Theta})}^{\circ}$\,(3.11). Now we will prove assertions (i) and (ii) together. The kernel $\osK$ of the homomorphism $\overline{\rho}_{F',F}: \osH_{F}^{\circ}\rightarrow \osH_{F'}^{\circ}$ is a smooth unipotent $\kappa$-subgroup, and the image $\osQ$   is a pseudo-parabolic $\kappa$-subgroup of  $\osH_{F'}^{\circ}$\,[P2,\,1.10\,(1),\,(2)].  The pseudo-parabolic subgroup $\osQ$ is clearly $\Theta$-stable as the facets $F$ and $F'$ are $\Theta$-stable. The kernel of  \,${\overline{\rho}}^{G}_{\cF',\cF}$ is $\osK\cap \osG^{\circ}_{\cF}$, and its image is contained in ${(\osQ^{\Theta})}^{\circ}$. Therefore, the kernel of ${\overline{\rho}}^{G}_{\cF',\cF}$ contains ${(\osK^{\Theta})}^{\circ}$ and is contained in $\osK^{\Theta}$.  As $\osK^{\Theta}$ is a smooth subgroup of $\osK$, we see that the kernel of  ${\overline{\rho}}^{G}_{\cF',\cF}$ is smooth.    
 
 \vskip1mm
 Since the image of the Lie algebra homomorphism ${\rm{L}}(\osG_{\cF}^{\circ})\rightarrow {\rm{L}}(\osG_{\cF'}^{\circ})$ induced by ${\overline\rho}^G_{\cF',\cF}$ is ${\rm{L}}(\osQ)^{\Theta}$, the containment $\fp(\cF'/\cF)= \overline{\rho}^G_{\cF',\cF}(\osG_{\cF}^{\circ}) \subset {(\osQ^{\Theta})}^{\circ}$ is equality. According to Lemma 3.7(ii), ${({\osQ}^{\Theta})}^{\circ}$ is a pseudo-parabolic $\kappa$-subgroup of $\osG^{\circ}_{\cF'}$. 
   \vskip1mm

     To prove  (iii), let $F'\prec F$ be as in the proof of (i) above and $\osQ$ be the image of \,$\overline{\rho}_{F',F}: \osH_F^{\circ}\rightarrow \osH_{F'}^{\circ}$. Then, as we saw above,  $\osQ$ is a $\Theta$-stable pseudo-parabolic $\kappa$-subgroup of $\osH^{\circ}_{F'}$ and $\fp(\cF'/\cF)=\osP := {(\osQ^{\Theta})}^{\circ}$. The inverse image of the subgroup $\osQ(\kappa)$  of $\osH^{\circ}_{F'}(\kappa)$ under the natural surjective homomorphism  $\sH^{\circ}_{F'} (\cO)\rightarrow \osH^{\circ}_{F'}(\kappa)$ equals  $\rho_{F',F}(\sH^{\circ}_{F}(\cO))\,(\subset \sH^{\circ}_{F'}(\cO))$, see [P2,\:1.10\,(4)].  Let $\sG_{F} ={(\sH^{\circ}_{F})}^{\Theta}$ and $\sG_{F'} ={(\sH^{\circ}_{F'})}^{\Theta}$. We will denote the $\cO$-group scheme homomorphism $\sG_{F}\rightarrow \sG_{F'}$  induced by $\rho_{F',F}$ by $\rho^{\Theta}_{F',F}$; the corresponding homomorphism $ \osG_{F}\rightarrow \osG_{F'}$ between the special fibers of \,$\sG_{F}$ and $\sG_{F'}$  will be denoted by $\overline{\rho}^{\Theta}_{F',F}$. The neutral components of $\sG_{F}$ and $\sG_{F'}$ are  $\sG^{\circ}_{\cF}$ and  $\sG^{\circ}_{\cF'}$ respectively (3.11).  Let $\sG^{\natural}_{F}\,(\supset \sG^{\circ}_{\cF})$ be the inverse image of $\sG^{\circ}_{\cF'}$ in $\sG_{F}$ under $\rho^{\Theta}_{F',F}$. Since  the homomorphism $\rho_{F',F}$ is the identity on the generic fiber $H$, we infer that $h\in \sH^{\circ}_{F}(\cO)$ is fixed under $\Theta$ if and only if so is $\rho_{F',F}(h)$, and as the generic fiber of both $\sG^{\circ}_{\cF}$ and $\sG^{\circ}_{\cF'}$ is $G$, the generic fiber of $\sG^{\natural}_{F}$ is also $G$.  It is easily seen now that the inverse image of the subgroup $\fp(\cF'/\cF)(\kappa)$ of \,$\osG^{\circ}_{\cF'}(\kappa)$, under the natural surjective homomorphism $\sG^{\circ}_{\cF'} (\cO)\rightarrow \osG^{\circ}_{\cF'}(\kappa)$, is  $\rho^{\Theta}_{F',F}(\sG^{\natural}_{F}(\cO))$.  We will presently show that the last group equals $\rho^G_{\cF',\cF}(\sG^{\circ}_{\cF}(\cO))$, this will prove (iii).    
\vskip1mm
 
  $\sG^{\natural}_{F}$ is the union of its generic fiber $G$ and its special fiber $\osG^{\natural}_{F}$; and the identity component of $\osG^{\natural}_{\cF}$ is clearly $\osG^{\circ}_{\cF}$.  We have shown above that the image $\osP$ of  \,$\osG^{\circ}_{\cF}$ under the homomorphism $\overline{\rho}^G_{\cF',\cF}$  is a pseudo-parabolic $\kappa$-subgroup of $\osG^{\circ}_{\cF'}$ and the kernel of this homomorphism is smooth. Hence, as $\kappa$ is separably closed, $\overline{\rho}^G_{\cF',\cF}(\osG^{\circ}_{\cF}(\kappa)) = \osP(\kappa)$.   So,  according to [CGP,\,Thm.\,C.2.23], there is a pseudo-parabolic $\kappa$-subgroup $\osP'$ of $\osG^{\circ}_{\cF'}$, that contains $\osP$, such that $\overline{\rho}^{\Theta}_{F',F}(\osG^{\natural}_{F}(\kappa)) = \osP'(\kappa)$. But since $\kappa$ is infinite, $\osP'(\kappa)/\osP(\kappa)$ is infinite unless $\osP' = \osP$. So we conclude that $\osP' = \osP$, and then $\overline{\rho}^{\Theta}_{F',F}(\osG^{\natural}_{F}(\kappa)) = \osP(\kappa)= \overline{\rho}^G_{\cF',\cF}(\osG^{\circ}_{\cF}(\kappa))$. Now using this, and the fact that the natural homomorphism $\sG^{\circ}_{\cF}(\cO)\rightarrow \osG^{\circ}_{\cF}(\kappa)$ is surjective (since $\cO$ is henselian and $\sG^{\circ}_{\cF}$ is smooth, [EGA\,IV$_4$, 18.5.17]) and the kernel of this homomorphism equals the kernel of the natural surjective homomorphism $\sG^{\natural}_{F}(\cO)\rightarrow \osG^{\natural}_{F}(\kappa)$, we  see that  $\rho^G_{\cF',\cF}(\sG^{\circ}_{\cF}(\cO))= \rho^{\Theta}_{F',F}(\sG^{\natural}_{F}(\cO))$. This proves (iii). 
 \vskip1mm 
  
   Finally, to prove the last assertion of the proposition, we fix a facet $F'$ of $\cB(H/K)$ that meets $\cF'$. Then $\sG^{\circ}_{\cF'}= \sG^{\circ}_{F'}$ (3.11).  Using Lemma 3.7(iii) for $\kappa$ in place of $k$ and $\osH^{\circ}_{F'}$ in place of $\cH$, we find a $\Theta$-stable pseudo-parabolic $\kappa$-subgroup $\osQ$ of $\osH^{\circ}_{F'}$ such that $\osP = {(\osQ^{\Theta})}^{\circ}$. Let $(F'\prec)\,F$ be the facet of $\cB(H/K)$ corresponding to the pseudo-parabolic $\kappa$-subgroup  $\osQ$ of $\osH^{\circ}_{F'}$. Then $F$ is stable under $\Theta$-action. As $F'\prec F$, there is a natural $\cO$-group scheme homomorphism $\rho_{F', F}: \sH^{\circ}_F\rightarrow \sH^{\circ}_{F'}$ that restricts to a $\cO$-group scheme homomorphism $\rho^G_{F', F}: \sG^{\circ}_F\rightarrow \sG^{\circ}_{F'}$.  Let $\osQ$ be the image of the former. Then according to (ii), the image of the latter is ${(\osQ^{\Theta}})^{\circ}=\osP$.  Let $P = \sG^{\circ}_F(\cO)\subset \sG^{\circ}_{F'}(\cO) =:Q$, and $\cF = \cF_P$. Then $P\subset Q$ are parahoric subgroups  of $G(K)$, $\cF'=\cF_Q\subset {\overline\cF}_Q\subset {\overline\cF}_P=\overline \cF$, thus $\cF'\prec \cF$. As $F$ and $F'$ meet $\cF$ and $\cF'$ respectively, $\sG^{\circ}_{\cF} = \sG^{\circ}_{F}$ and $\sG^{\circ}_{\cF'} = \sG^{\circ}_{F'}$\,(3.11), and hence the image of the homomorphism $\overline{\rho}^G_{\cF',\cF}: \osG^{\circ}_{\cF}\rightarrow   \osG^{\circ}_{\cF'}$ equals $\osP$. \hfill$\Box$   
 
 \vskip2mm

Proposition 3.12 and [CGP, Propositions 2.2.10 and 3.5.1] imply the following. (Recall that the residue field $\kappa$ of $K$ has been assumed to be separably closed!)
 
\vskip2mm
\ni{\bf Corollary 3.13.} {\em {\rm(i)}  A facet $\cF$ of $\sB$ is a chamber\,$($=maximal facet$)$ if and only if \,$\osG^{\circ}_{\cF}$ does not contain a proper pseudo-parabolic $\kappa$-subgroup. Equivalently, $\cF$ is a chamber if and only if the pseudo-reductive quotient \,$\orG_{\cF}:= \osG^{\circ}_{\cF}/\sR_{u,\kappa}(\osG^{\circ}_{\cF})$ of \,$\osG^{\circ}_{\cF}$  is commutative $($this is the case if and only if \,$\orG_{\cF}$  contains a unique maximal $\kappa$-torus, or, equivalently, every torus of this pseudo-reductive group is central\,$)$. 
\vskip1mm
{\rm(ii)}The codimension of a facet $\cF$ of $\sB$ equals the rank of the derived subgroup of  the pseudo-split  pseudo-reductive quotient \,$\orG_{\cF}$. Hence the dimension of $\cF$ equals the dimension of the maximal central torus of $\orG_{\cF}$\,$($recall that the central torus of $G$ is $K$-anisotropic\,{\rm (2.13)},  the dimension of $\sB$ equals $K$-{\rm rank}\:$G$\,{\rm (2.16)}, and {\rm rank}\,$\orG_{\cF}$= $K$-{\rm rank}\:$G$$)$. Thus a point $x$ of $\sB$ is a vertex $($of $\sB )$ if and only if $\orG_x := \osG^{\circ}_x/\sR_{u,\kappa}(\osG^{\circ}_x)$ does not contain a nontrivial central torus.} 
\vskip2mm

We  will now establish the following analogues of Propositions 3.5--3.7 of [P2].  
\vskip2mm

\ni{\bf Proposition 3.14.} {\em Let $\cA$ be an apartment  of $\sB$, and $\cC$, $\cC'$ two chambers in $\cA$. Then there is a {gallery} joining $\cC$ and $\cC'$ in $\cA$, i.e., there is a finite sequence 
$$\cC=\cC_0, \:\cC_1, \:\ldots\,,\:\cC_m =\cC'$$ 
of chambers in $\cA$ such that for $i$ with $1\leqslant i\leqslant m$, $\cC_{i-1}$ and $\cC_i$ share a face of codimension 1.}
\vskip2mm

\ni{\it Proof.}  Let $\cA_2$ be the codimension $2$-skelton of $\cA$, i.e.,  the union of all facets in $\cA$ of codimension at least $2$. Then $\cA_2$ is a closed subset of $\cA$ of codimension $2$, so $\cA-\cA_2$ is a connected open subset of the affine space $\cA$. Hence $\cA-\cA_2$ is arcwise connected. This implies that given points $x\in \cC$ and $x'\in \cC'$, there is a piecewise linear curve in $\cA-\cA_2$ joining $x$ and $x'$. Now the chambers in $\cA$ that meet this curve make a gallery joining $\cC$ to $\cC'$. \hfill$\Box$
\vskip1mm

As the central torus of $G$ is $K$-anisotropic, the dimension of any apartment, or any chamber, in $\sB$ is equal to the $K$-rank of $G$. 
A {\it panel} in $\sB$ is by definition a facet of codimension $1$.
\vskip2mm

\ni{\bf Proposition 3.15.} {\em $\sB$ is thick, that is any panel is a face of at least three chambers, 
and every apartment of $\sB$ is thin, that is any panel lying in an apartment is a face of exactly two chambers of the apartment.}
\vskip2mm

\ni{\it Proof.} Let $\cF$ be a  facet of $\sB$ that is not a chamber, and $\cC$ be a chamber of which $\cF$ 
is a face. Then there is an $\cO$-group scheme homomorphism $\rho^G_{\cF,\cC}: \sG^{\circ}_{\cC}\rightarrow \sG^{\circ}_{\cF}$ (3.2). The image of $\osG^{\circ}_{\cC}$ in $\osG^{\circ}_{\cF}$,  under the induced homomorphism of special fibers, is a minimal pseudo-parabolic $\kappa$-subgroup of $\osG^{\circ}_{\cF}$, and conversely, any minimal pseudo-parabolic $\kappa$-subgroup of the latter determines a chamber of $\sB$ with $\cF$ as a face (Corollary 3.13). Now as $\kappa$ is infinite, $\osG^{\circ}_{\cF}$ contains infinitely many minimal pseudo-parabolic $\kappa$-subgroups.  We conclude that $\cF$ is a face of infinitely many chambers.
\vskip1mm

%
The second assertion follows at once from the following well-known result in algebraic topology:\:In any simplicial complex whose geometric realization is a  topological manifold without boundary (such as an apartment $\cA$ of $\sB$), any simplex of codimension 1 is a face of exactly two chambers (i.e., maximal dimensional simplices).\hfill$\Box$  

\vskip2mm

\ni{\bf Proposition 3.16.} {\em Let $\cA$ be an apartment of  $\sB$ and $S$ be the maximal $K$-split torus of $G$ corresponding to this apartment. $($Then $\cA = \cB(Z_H(S)/K)^{\Theta}.)$ The group $N_G(S)(K)$ acts transitively on the set of chambers of $\cA$.} 
\vskip2mm

\ni{\it Proof.} According to  Proposition 3.14, given any  two chambers in $\cA$, there exists a minimal gallery in $\cA$ joining these two chambers. So to prove the proposition by induction on the length of a minimal gallery joining two chambers, it suffices to prove that given two different chambers  $\cC$ and $\cC'$ in $\cA$ which share a panel $\cF$, there is an element $n\in N_G(S)(K)$ such that $n\cdot \cC = \cC'$. Let $\sG:=\sG^{\circ}_{\cF}$ be the Bruhat-Tits smooth affine $\cO$-group scheme associated to the panel $\cF$ and $\sS\subset \sG$ be the closed $\cO$-torus with generic fiber $S$. Let  $\osG$ be the special fiber of $\sG$ and $\osS$ the special fiber of $\sS$. Then  $\osS$ is a maximal torus of $\osG$. The chambers $\cC$ and $\cC'$ correspond to minimal pseudo-parabolic subgroups $\osP$ and $\osP'$ of $\osG$ (Corollary 3.13). Both of  these minimal pseudo-parabolic $\kappa$-subgroups contain $\osS$ since the chambers $\cC$ and $\cC'$ lie in $\cA$. But then by Theorems C.2.5 and C.2.3 of [CGP], there is an element $\overline{n}\in \osG(\kappa)$ that normalizes $\osS$ and conjugates $\osP$ onto $\osP'$.   Now from Proposition 2.1(iii) of [P2] we conclude that there is an element $n\in N_{\sG}(\sS)(\cO)$ lying over $\overline{n}$.  It is clear that $n$ normalizes $S$ and hence it lies in $N_G(S)(K)$; it fixes $\cF$ pointwise and $n\cdot\cC =\cC'$.\hfill$\Box$   
\vskip2mm

Now in view of Propositions 2.17, 3.4, 3.14 and 3.15, Theorem 3.11 of [Ro] (cf.\,also [P2,\,1.8]) implies that  $\sB$ is an affine building if  for any maximal $K$-split torus  $S$ of $G$, $\cB(Z_H(S)/K)^{\Theta}$ is taken to be the corresponding apartment, and  $\sB$ is given the   polysimplicial structure described in 3.5. Thus we obtain the following:
\vskip2mm

\ni{\bf Theorem 3.17.} {\em $\sB=\cB(H/K)^{\Theta}$ is an affine building. Its apartments are the affine spaces 
$\cB(Z_H(S)/K)^{\Theta}$ under $V(S) :=\bR\otimes_{\bZ}{\rm{X}}_{\ast}(S)$ for maximal $K$-split tori $S$ of $G$.  Its facets are as in $3.5$. 
The group $G(K)$ acts on $\sB$ by polysimplicial isometries.}   
\vskip3mm

From Propositions 2.18 and 3.16 we obtain the following.
\vskip2mm

\ni{\bf Proposition 3.18.} {\em $G(K)$ acts transitively on the set of ordered pairs $(\cA,\cC)$ consisting of an apartment $\cA$ of $\sB$ and a chamber $\cC$ of $\cA$.} 
\vskip3mm

\ni{\bf Remark 3.19.}  (i) As in [P2,\,3.16], using the preceding proposition we can obtain Tits systems in suitable subgroups of $G(K)$. 
\vskip1mm

(ii) As in [P2, \S5], we can obtain filtration of root groups and a valuation of root datum for $G/K$. 

\vskip5mm
\ni{\bf \S4. Tamely-ramified descent} 

\vskip3mm
We begin by proving the following proposition in which $\kappa$ is any field of characteristic $p\geqslant 0$. 
\vskip2mm

\ni{\bf Proposition 4.1.} {\em Let $\cH$ be a {noncommutative} pseudo-reductive $\kappa$-group, $\theta$ a $\kappa$-automorphism of $\cH$ of finite order not divisible by $p$, and $\cG := {(\cH^{\langle \theta\rangle})}^{\circ}$. Then} 
\vskip1mm

{\rm{(i)} {\em No maximal torus of $\cG$ is central in $\cH$.  } 
\vskip1mm

{\rm{(ii)} {\em The centralizer in $\cH$ of any maximal torus of $\cG$ is commutative.} 
\vskip1mm

{\rm{(iii)}} {\em Given a maximal $\kappa$-torus $\cS$ of $\cG$, there is a $\theta$-stable maximal $\kappa$-torus of $\cH$ containing $\cS$.} 
\vskip1mm

{\rm{(iv)} {\em If $\kappa$ is separably closed, then $\cH$ contains a $\theta$-stable proper pseudo-parabolic $\kappa$-subgroup.}
\vskip2mm

\ni{\it Proof.}  We fix an algebraic closure $\overline\kappa$ of $\kappa$.  Let $\cH'$ be the maximal reductive quotient of $\cH_{\overline\kappa}$. As $\cH$ is noncommutative, $\cH'$ is also noncommutative (see [CGP, Prop.\,1.2.3]). The automorphism $\theta$ induces a $\overline\kappa$-automorphism of $\cH'$ which we will denote again by $\theta$. 
According to a theorem of Steinberg\,[St, Thm.\,7.5], $\cH_{\overline\kappa}$ contains a $\theta$-stable Borel subgroup $\sB$, and this Borel subgroup contains a $\theta$-stable maximal torus $\cT$. The natural quotient map $\pi: \cH_{\overline\kappa}\rightarrow \cH'$ carries $\cT$ isomorphically onto a maximal torus of $\cH'$.  We endow the root system of $\cH'$ with respect to the maximal torus $\fT':= \pi(\cT)\cap \sD(\cH')$  of the derived subgroup $\sD(\cH')$ of $\cH'$ with the ordering determined by the Borel subgroup $\pi(\cB)$. Let $a$ be the sum of all positive roots. Then as $\pi(\cB)$ is $\theta$-stable, $a$ is fixed under $\theta$ acting on the character group ${\rm X}(\fT')$ of $\fT'$. Therefore, ${\rm X}(\fT')$ admits a nontrivial torsion-free quotient on which $\theta$ acts trivially. This implies that $\cT$ contains a nontrivial subtorus $\sT$  that is fixed pointwise under $\theta$ and is mapped by $\pi$ into $\fT'\,(\subset \sD(\cH'))$. The subtorus $\sT$  is therefore contained in $\cG_{\overline\kappa}$.  Since the center of the semi-simple group $\sD(\cH')$ does not contain a nontrivial smooth connected subgroup, we infer that $\sT$ is not central in $\cH_{\overline\kappa}$.  Thus the subgroup $\cG_{\overline\kappa}$ contains a  noncentral  torus of $\cH_{\overline\kappa}$. Now by conjugacy of maximal tori in $\cG_{\overline\kappa}$, we see that no maximal torus of this group can be central in $\cH_{\overline\kappa}$. This proves (i).  
\vskip.5mm

To prove (ii), it would suffice to show that the centralizer $Z_{\cH}(\cS)$ in $\cH$ of a maximal $\kappa$-torus $\cS$ of $\cG$  is commutative.  This centralizer is a $\theta$-stable pseudo-reductive subgroup of $\cH$, and ${(Z_{\cH}(\cS)^{\langle \theta\rangle})}^{\circ} = Z_{\cG}(\cS)$. As $\cS$ is a maximal torus of  $Z_{\cG}(\cS)$ that is central in $Z_{\cH}(\cS)$, if  $Z_{\cH}(\cS)$ were noncommutative, we could apply (i) to this subgroup in place of $\cH$ to get a contradiction.     
\vskip.5mm

To prove (iii), we consider the centralizer $Z_{\cH}(\cS)$ of $\cS$ in $\cH$. This centralizer  is $\theta$-stable and commutative according to\,(ii). The unique maximal $\kappa$-torus of  it contains $\cS$ and is a $\theta$-stable maximal torus of $\cH$.
\vskip.5mm

To prove (iv), we assume now that $\kappa$ is separably closed and let $\cS$ be a maximal torus of $\cG$. Then $\cS$ is $\kappa$-split, and in view of (i), there is a 1-parameter subgroup $\lambda:{\rm{GL}}_1\rightarrow \cS$ whose image is not central in $\cH$.  Then $P_{\cH}({\lambda})$ is a $\theta$-stable proper pseudo-parabolic $\kappa$-subgroup of $\cH$. \hfill$\Box$  
    
 \vskip2mm
 
 In the following proposition we will use the notation introduced in \S\S1,\,2. As in 2.13, we will assume that $H$ is semi-simple and the central torus of $G$ is $K$-anisotropic.  We will further assume that $H$ is $K$-isotropic, $\Theta$ is a finite {\it cyclic} group of automorphisms of $H$, and $p$ does not divide the order of $\Theta$.    
 
 \vskip2mm
  
\ni{\bf Proposition 4.2.}  {\em The Bruhat-Tits building $\cB(H/K)$ of $H(K)$ contains a $\Theta$-stable chamber.} 
\vskip2mm

\ni{\it Proof.}  Let  $F$ be a $\Theta$-stable facet of $\cB(H/K)$ that is maximal among the $\Theta$-stable facets. Let $\sH := \sH^{\circ}_F$ be the Bruhat-Tits smooth affine $\cO$-group scheme with generic fiber $H$, and connected special fiber $\osH$, corresponding to $F$.  Let  $\cH:= \osH/\sR_{u,\kappa}(\osH)$ be the maximal pseudo-reductive quotient of $\osH$. In case $\cH$ is commutative, $\osH$ does not contain a proper pseudo-parabolic $\kappa$-subgroup and so $F$ is a chamber of $\cB(H/K)$. We assume, if possible, that $\cH$ is not commutative.  As $F$ is stable under the action of $\Theta$, there is a natural action of this finite cyclic group on $\sH$ by $\cO$-group scheme automorphisms\,(2.5(i)). This action induces an action of $\Theta$ on $\osH$, and so also on its pseudo-reductive quotient $\cH$. Now taking $\theta$ to be a generator of $\Theta$, and using Proposition 4.1(iv), we conclude that $\cH$ contains a $\Theta$-stable proper pseudo-parabolic $\kappa$-subgroup. The inverse image $\osP$ in $\osH$ of any such pseudo-parabolic subgroup of $\cH$ is a $\Theta$-stable proper  pseudo-parabolic $\kappa$-subgroup of $\osH$.   The facet $F'$ corresponding to $\osP$ is $\Theta$-stable and $F\prec F'$. This contradicts the maximality of $F$. Hence, $\cH$ is commutative and $F$ is a chamber.\hfill$\Box$  

\vskip2mm

To prove the next theorem\,(Theorem 4.4), we will use the following:
\vskip2mm

\ni{\bf Proposition 4.3.} {\em Let $\fK$ be a field complete with respect to a discrete valuation and with separably closed residue field. Let $\fG$ be a connected absolutely simple $\fK$-group of inner type $A$ that splits over a finite tamely-ramified field extension $\fL$ of $\fK$. Then $\fG$ is $\fK$-split.}
\vskip2mm 

\ni{\it Proof.} We may (and do) assume that $\fG$ is simply connected.  Then $\fG$ is $\fK$-isomorphic to ${\rm{SL}}_{n,\fD}$, where $\fD$ is a finite dimensional division algebra with center $\fK$ that splits over the  finite tamely-ramified field extension $\fL$ of $\fK$. By Propositions 4 and 12 of [S, Ch.\,II] the degree of $\fD$ is a power of $p$, where $p$ is the characteristic of the residue field of $\fK$. But a noncommutative division algebra of degree a power of $p$ cannot split over a field extension of degree prime to $p$. So, $\fD =\fK$, hence $\fG \simeq {\rm{SL}}_n$ is $\fK$-split. \hfill$\Box$   
\vskip2mm

\ni{\bf Theorem 4.4.} {\em A  semi-simple $K$-group $G$  that is quasi-split over a finite tamely-ramified field extension of $K$ is already quasi-split over $K$.}
\vskip1mm
  
  This theorem has been proved by Philippe Gille in [Gi] by an entirely different method. 
\vskip2mm

\ni{\it Proof.} We assume that all field extensions appearing in this proof are contained in a fixed separable closure of $K$. To prove the theorem, we may (and do) replace $G$ by its simply-connected central cover and assume that $G$ is simply connected. Let $S$ be a maximal $K$-split torus of $G$. Then $G$ is quasi-split over a (separable)  extension $L$ of $K$ if and only if the derived subgroup $Z_G(S)'$ of the centralizer $Z_G(S)$ of  $S$ is quasi-split over $L$. Moreover, $G$ is quasi-split over $K$ if and only if $Z_G(S)'$ is trivial. Therefore, to prove the theorem we need to show that any semi-simple simply connected $K$-{\it anisotropic} $K$-group that is quasi-split over a finite tamely-ramified field extension of $K$ is necessarily trivial.   Let $G$ be any such group. In view of Proposition 1.1, we may replace $K$ by its completion and  assume that $K$ is complete. 
\vskip1mm

There exists a finite indexing set $I$, and for each $i\in I$,  a finite  separable field extension $K_i$ of $K$ and an absolutely almost simple simply connected $K_i$-anisotropic $K_i$-group $G_i$  such that $G = \prod_{i\in I}{\rm{R}}_{K_i/K}(G_i)$. Now $G$ is quasi-split over a finite separable field extension $L$ of $K$ if and only if for each $i$,  ${\rm{R}}_{K_i/K}(G_i)$ is quasi-split over $L$.  But ${\rm{R}}_{K_i/K}(G_i)$ is quasi-split over $L$ if and only if $G_i$ is quasi-split over the compositum $L_i := K_iL$. For  $i\in I$,  the finite extension $K_i$ of $K$ is complete and its residue field is separably closed, and if $L$ is a finite tamely-ramified field extension of $K$, then  $L_i$ is a finite tamely-ramified field extension of $K_i$. So to prove the theorem, we may (and do) replace $K$ by $K_i$ and $G$ by $G_i$ to assume that $G$ is an absolutely almost simple simply connected {\it $K$-anisotropic} $K$-group that is quasi-split over a finite tamely-ramified field extension of $K$. We will show that such a group $G$  is trivial.  
\vskip1mm

Since the residue field $\kappa$ of $K$ is separably closed, any finite tamely-ramified field extension of $K$ is a cyclic Galois extension. Let $L$ be a finite tamely-ramified field extension of $K$ of minimal degree over which $G$ is quasi-split.  The Galois group $\Theta$ of $L/K$ is a finite cyclic group and its order is not divisible by $p$\,(= char($\kappa$)).  
 \vskip1mm
 
 As $G_L$ is quasi-split, Bruhat-Tits theory is available for $G$ over $L$ [BrT2,\,\S4].   Let $H = {\rm{R}}_{L/K}(G_L)$. Then $H$ is quasi-split over $K$ and hence Bruhat-Tits theory is also available for $H$ over $K$. Let $\cB(H/K)$ be the Bruhat-Tits building of $H(K)\,(=G(L))$. Elements of  $\Theta$ act by $K$-automorphisms on $H$ and so on $\cB(H/K)$ by polysimplicial isometries; moreover, $G = H^{\Theta}$. The Galois group $\Theta$ acts on $G(L)$ by continuous automorphisms and so it acts on the Bruhat-Tits building $\cB(G/L)$ of $G(L)$ by polysimplicial isometries.  There is a natural $\Theta$-equivariant identification of the building $\cB(H/K)$ with the building $\cB(G/L)$. (Note that $K$-rank\,$H$ = $L$-rank\,$G_L$, and there is a natural bijective correspondence between the set of maximal $K$-split tori of $H$ and the set of maximal $L$-split tori of $G_L$, see [CGP, Prop.\,A.5.15(2)]. This correspondence will be used below.)  The results of \S3 imply that Bruhat-Tits theory is available for $G$ over $K$ and $\sB:=\cB(H/K)^{\Theta} \,(= \cB(G/L)^{\Theta})$ is the Bruhat-Tits building of $G(K)$. 
 \vskip1mm
 
 Since $G$ is $K$-anisotropic, the building $\sB$ of $G(K)$ consists of a single point (Proposition 2.17), hence $\Theta$ fixes a unique point of  \,$\cB(G/L)$. Let $C$ be the facet of $\cB(G/L)$ that contains this point.  Then $C$ is stable under $\Theta$. According to Proposition 4.2, $C$ is a chamber. Let $\sH := \sH_{C}^{\circ}$ be the Bruhat-Tits smooth affine $\cO$-group scheme associated to $C$ with generic fiber $H$ and connected special fiber $\osH$.  As $C$ is a chamber,  the maximal pseudo-reductive quotient $\osH^{\rm{pred}}$ of $\osH$  is commutative [P2,\:1.10].  Now using Proposition 2.6 for $\Omega = C$ we obtain a $\Theta$-stable maximal $K$-split torus $T$ of $H$ such that $C$ lies in the apartment $A(T)$ corresponding to $T$\:(and the special fiber of the schematic closure of $T$ in $\sH$ maps onto the maximal torus of $\osH^{\rm{pred}}$).  Let $T'$ be the image of $T_L$ under the natural surjective homomorphism $q:\,H_L = {\rm{R}}_{L/K}(G_L)_L\rightarrow G_L$\,(for the proof of surjectivity of $q$, see [CGP, Prop.\,A.5.11(1)]). Then  $T'$ is a $L$-torus of $G_L$ and according to [CGP,\:Prop.\,A.5.15(2)] it is the unique maximal $L$-split torus of $G_L$ such that ${\rm{R}}_{L/K}(T')\,(\subset {\rm{R}}_{L/K}(G_L)=H)$ contains the maximal $K$-split torus  $T$ of $H$. 
 \vskip1mm
 
 We identify $H(K)$ with $G(L)$. Then for $x\in H(K) (\subset H(L))$ and $\theta\in \Theta$, we have $q(\theta(x)) = \theta(x)$. Since $T(K)$ is $\Theta$-stable,  for $t\in T(K)$ and $\theta\in \Theta$, $\theta(t)$ lies in $T'(L)$. Now as $T(K)$ is Zariski-dense in $T$, its image in $T'(L)$ is Zariski-dense in $T'$. Since this  image is stable under the action of $\Theta = {\rm{Gal}}(L/K)$ on $G(L)$, from the Galois criterion\,[Bo, Ch.\,AG, Thm.\,14.4(3)] we infer that $T'$ descends to a $K$-torus  of $G$, i.e., there is a $K$-torus $\cT$ of $G$ such that $T' = \cT_L$.    
 In the natural identification of $\cB(H/K)$ with $\cB(G/L)$, the apartment $A(T)$ of the former is $\Theta$-equivariantly  identified with the apartment $A(T')$ of the latter. We will view the chamber $C$ as a $\Theta$-stable chamber in $A(T')$.  
 \vskip1mm
 
 Let $\Delta$ be the basis of the affine root system of the absolutely almost simple, simply connected quasi-split  $L$-group $G_L$ with respect to $T' \,(= \cT_L)$, determined by the $\Theta$-stable chamber $C$ [BrT2, \S4]. Then  $\Delta$ is stable under the action of $\Theta$ on the affine root system of $G_L$ with respect to $T'$. There is a natural  $\Theta$-equivariant bijective correspondence between the set of vertices of $C$ and $\Delta$. Since $\sB$, and hence $C^{\Theta}$,  consists of a single point, $\Theta$ acts transitively on the set of vertices of $C$ so it acts transitively on $\Delta$. Now from the classification of irreducible affine root systems [BrT1,\,\S1.4.6], we see that $G_L$ is a {\it split} group of type $A_n$ for some $n$.  Proposition 4.3 implies that $G$ cannot be of inner type $A_n$ over $K$. On the other hand, if $G$ is of outer type $A_n$, then over a quadratic Galois extension $K'\,(\subset L)$ of $K$ it is of inner type.  Now, according to Proposition 4.3, $G$  splits over $K'$. We conclude that $L=K' $ and hence $\#\Theta = 2$. As $\Theta$ acts transitively on $\Delta$ and $\#\Delta = n+1$, we infer that $n+1 = 2$, i.e., $n=1$, and then $G$ is of inner type, a contradiction. \hfill$\Box$   
 \vskip1mm

\vskip2mm
 
\ni{\bf 4.5.} Now let $k$ be a field endowed with a nonarchimedean discrete valuation. We assume that  the valuation ring of $k$ is Henselian. Let $K$ be the maximal unramified extension of $k$, and $L$ be a finite tamely-ramified field extension of $K$ with Galois group $\Theta :={\rm{Gal}}(L/K)$. Let $G$ be a connected reductive $k$-group that is quasi-split over $K$ and $H = {\rm{R}}_{L/K}(G_L)$. Then $G = H^{\Theta}$, and by Theorem 3.17, the Bruhat-Tits building $\cB(G/K)$ of $G(K)$ can be identified with the subspace of points in the Bruhat-Tits building of $G(L) \,(= H(K))$ that are fixed under $\Theta$\,(with polysimplicial structure on $\cB(G/K)$ as in 3.5).  Now by ``unramified descent''\,[P2], {\em Bruhat-Tits theory is available for $G$ over $k$ and the Bruhat-Tits building of $G(k)$ is $\cB(G/K)^{{\rm{Gal}}(K/k)}$.} 
\vskip4mm

 \centerline{\bf References}
 
\vskip2mm
 

 \ni[Bo] A.\,Borel, {\it Linear algebraic groups} (second edition). Springer-Verlag, New York (1991).
 \vskip1mm
 
 
 
 \ni[BLR] S.\,Bosch, W.\,L\"utkebohmert and M.\,Raynaud, {\it N\'eron models}. Springer-Verlag, Heidelberg (1990).
\vskip1mm
 
 \ni[BrT1] F.\,Bruhat and J.\,Tits, {\it Groupes r\'eductifs sur un corps local.} Publ.\,Math.\,IHES {\bf 41}(1972).
 \vskip1mm
 
 \ni[BrT2] F.\,Bruhat and J.\,Tits, {\it Groupes r\'eductifs sur un corps local, II.} Publ.\,Math.\,IHES {\bf 60}(1984).
 \vskip1mm
 
   

%

 \ni[CGP] B.\,Conrad, O.\,Gabber and G.\,Prasad, {\it Pseudo-reductive groups} (second edition). Cambridge U.\,Press, New York (2015).
 \vskip1mm 
 
\ni[$\rm{SGA3_{II}}$] M.\,Demazure and A.\,Grothendieck, {\it Sch\'emas en groupes, Tome II}. Lecture Notes in Math.\,\#{\bf 152}, Springer-Verlag, Heidelberg (1970). 
 \vskip1mm
 

\ni[E] B.\,Edixhoven, {\it N\'eron models and tame ramification}. Comp.\,Math.\,{\bf 81}(1992), 291-306. 
\vskip1mm

\ni[GGM] O.\,Gabber, P.\,Gille and L.\,Moret-Bailly, {\it Fibr\'es principaux sur les corps valu\'es hens\'eliens.} Algebr.\,Geom.\,{\bf 1}(2014), 573-612. 
\vskip1mm
%
%

\ni[Gi] P.\,Gille, {\it Semi-simple groups that are quasi-split over a tamely-ramified extension}, Rendiconti Sem.\,Mat.\,Padova, to appear. 
\vskip1mm
 
  \ni[$\rm{EGA\,IV_3}$] A.\,Grothendieck, {\it {\it El\'ements de g\'eom\'etrie alg\'ebrique}, IV: \'Etude locale des sch\'emas et des morphismes de sch\'emas},  Publ.\,Math.\,IHES {\bf 28}(1966),  5-255. 
 \vskip1mm 

 \ni[$\rm{EGA\,IV_4}$] A.\,Grothendieck, {\it {\it El\'ements de g\'eom\'etrie alg\'ebrique}, IV: \'Etude locale des sch\'emas et des morphismes de sch\'emas},  Publ.\,Math.\,IHES {\bf 32}(1967),  5-361. 
 \vskip1mm
 
\ni[P1] G.\,Prasad, {\it Galois-fixed points in the Bruhat-Tits building of a reductive group}. Bull.\,Soc.\,Math.\,France\,{\bf 129}(2001), 169-174.
\vskip1mm

\ni[P2] G.\,Prasad, {\it A new approach to unramified descent in Bruhat-Tits theory}. American J.\,Math.\,(to appear).  
\vskip1mm
 
 \ni[PY1] G.\,Prasad and J.-K.\,Yu, {\it On finite group actions on reductive groups and buildings.} Invent.\,Math.\,{\bf 147}(2002), 545--560.
\vskip1mm

\ni[PY2] G.\,Prasad and J.-K.\,Yu, {\it On quasi-reductive group schemes.} J.\,Alg.\,Geom. 
{\bf 15}\,(2006), 507-549.
\vskip1mm 

\ni[Ri] R.\,Richardson, {\it On orbits of algebraic groups and Lie groups}, Bull.\,Australian Math.\,Soc.\,{\bf 25}(1982), 1-28.
\vskip1mm 

 \ni[Ro] M.\,Ronan, {\it Lectures on buildings}. University of Chicago Press, Chicago (2009).
 \vskip1mm
 
 \ni[Rou] G.\,Rousseau, {\it Immeubles des groupes r\'eductifs sur les corps locaux}, University of Paris, Orsay, thesis (1977).
 
 (Available at http://www.iecl.univ-lorraine.fr/{\raise.17ex\hbox{$\scriptstyle\sim$}}Guy.Rousseau/Textes/) 
\vskip1mm
  
\ni[S] J-P. Serre,  {\it Galois cohomology.} Springer-Verlag, New York (1997).
\vskip1mm

\ni[St] R.\,Steinberg, {\it Endomorphisms of linear algebraic groups}, Memoirs of the Amer.\,Math. Soc.\:{\bf 80}(1968).
\vskip1mm
 

\ni[T] J.\,Tits, {\it Reductive groups over local fields.} Proc.\,Symp.\,Pure 
Math.\,\#{\bf 33}, Part I, 29--69, American Math.\,Soc.\,(1979).
\vskip1mm


\vskip1mm

\ni{University of Michigan}

\ni{Ann Arbor, MI 48109.}

\ni{e-mail: gprasad@umich.edu}

\end{document}